\documentclass{amsart}
\usepackage{amsmath,amssymb,amscd,latexsym,amsthm,mathrsfs,comment}
\usepackage[unicode]{hyperref}
\newtheorem{theorem}{Theorem}

\newtheorem{Theorem}{Theorem}
\newtheorem{ATheorem}{Theorem}

\newtheorem{Lemma}{Lemma}
\newtheorem{ALemma}{Lemma}

\newtheorem*{corollary}{Corollary}

\theoremstyle{remark}
\newtheorem*{algo}{\bf Algorithm}
\newtheorem{Remark}{\bf Remark}
\newtheorem{example}{\bf Example}
\newtheorem*{demo}{Proof}
\let\rom\textup

\def\H{\mathbb H}
\def\P{\mathbb P}

\def\Z{\mathbb Z}
\def\R{\mathbb R}
\def\C{\mathbb C}
\def\CC{\mathscr C}
\def\adj{\mathrm{adj}}
\def\Re{\operatorname{Re}}
\def\Im{\operatorname{Im}}
\def\gcd{\operatorname{gcd}}
\def\Sp{\textup{Sp}}
\def\SL{\textup{SL}}
\def\GL{\textup{GL}}
\def\MF#1#2#3#4#5{\begin{pmatrix} \MMF #1 \\ \MMF #2 \\ \MMF #3 \\
  \MMF #4 \\ \MMF #5\end{pmatrix}}
\def\MMF#1#2#3#4#5{#1 & #2 & #3 & #4 & #5}

\def\M#1#2#3#4{\begin{pmatrix} \MM #1 \\ \MM #2 \\ \MM #3 \\
  \MM #4 \end{pmatrix}}
\def\MM#1#2#3#4{#1 & #2 & #3 & #4}
\let\<\langle
\let\>\rangle
\newcommand{\Tau}{{\mathrm T}}
\newcommand{\Zeta}{{\mathrm Z}}
\renewcommand{\d}{{\mathrm d}}

\newcommand{\wt}{\widetilde}
\newcommand{\wh}{\widehat}
\newcommand\trans[1]{{}^t#1}
\DeclareMathAlphabet{\mathpzc}{U}{euf}{m}{n}
\def\fa{\mathpzc a}
\def\fb{\mathpzc b}
\def\fc{\mathpzc c}
\def\fd{\mathpzc d}
\def\bw{\mathbf w}
\def\ba{\mathbf a}

\begin{document}
\hypersetup{pdfauthor={Yifan Yang and Wadim Zudilin},%
pdftitle={On Sp\_4 modularity of Picard--Fuchs differential equations for Calabi--Yau threefolds}}
\title[$\Sp_4$ modularity of Picard--Fuchs differential equations]
  {On $\Sp_4$ modularity of Picard--Fuchs differential equations for
    Calabi--Yau threefolds}
\author{Yifan Yang}
\address{Department of Applied Mathematics, National Chiao Tung
  University, Hsinchu 300, TAIWAN}
\email{yfyang@math.nctu.edu.tw}
\author{Wadim Zudilin}
\address{Max-Planck-Institut f\"ur Mathematik, Vivatsgasse 7, 53111 Bonn, GERMANY}
\email{wzudilin@mpim-bonn.mpg.de}
\date{October 10, 2008}

\thanks{The work of both authors was
  supported by fellowships of the Max Planck Institute for
  Mathematics (Bonn). The work of the first author was also supported by grant
  95-2115-M-009-005 of the National Science Council (NSC) of the
  Republic of China (Taiwan).}

\begin{abstract}
Motivated by the relationship of classical modular functions and
Picard--Fuchs linear differential equations of order 2 and~3, we present
an analogous concept for equations of order 4 and~5.
\end{abstract}

\maketitle

\section{Introduction}
\label{s0}

Let $M_z$ be a family of Calabi--Yau threefolds parameterized by a
complex variable $z\in\P^1(\C)$. Then periods of the unique
holomorphic differential $3$-form on~$M_z$ satisfy a linear
differential equation, called the \emph{Picard--Fuchs differential
equation} of~$M_z$. When the Hodge number $h^{2,1}$ is equal to~$1$,
the Picard--Fuchs differential equation has order~$4$. One of the most
well-known examples is perhaps the family of quintic threefolds
$$
  x_1^5+x_2^5+x_3^5+x_4^5+x_5^5-z^{-1/5}x_1x_2x_3x_4x_5=0
$$
in $\P^4$, whose Picard--Fuchs differential equation is
\begin{equation}
\label{eq: quintic}
  \theta^4y-5z(5\theta+1)(5\theta+2)(5\theta+3)(5\theta+4)y=0, \qquad
  \theta=z\frac\d{\d z}
\end{equation}
(see~\cite{COGP}). This is one of the fourteen families of Calabi--Yau
threefolds whose Picard--Fuchs differential equations are
hypergeometric (we refer the reader to the classical book~\cite{Sl}
for the definition of hypergeometric functions and hypergeometric
differential equations).

Very recently, we \cite{CYY} studied the monodromy aspect of Picard--Fuchs
differential equations originated from Calabi--Yau
threefolds. One of the main results in~\cite{CYY} is that,
with respect to certain bases, the monodromy groups of the fourteen
hypergeometric Picard--Fuchs differential equations are contained in
certain congruence subgroups of $\Sp_4(\Z)$ (see
\cite[Theorem~2]{CYY}). For instance, in the case of~\eqref{eq: quintic}, the
monodromy matrices around the singular points $z=0$ and $z=1/3125$ are
\begin{equation}
\label{eq: gen_quintic}
  \begin{pmatrix} 1 & 1 & 0 & 0 \\ 0 & 1 & 0 & 0 \\
  5 & 5 & 1 & 0 \\ 0 & -5 & -1 & 1 \end{pmatrix}
\quad\text{and}\quad
  \begin{pmatrix} 1 & 0 & 0 & 0 \\ 0 & 1 & 0 & 1 \\
  0 & 0 & 1 & 0 \\ 0 & 0 & 0 & 1\end{pmatrix},
\end{equation}
respectively. The group generated by these two matrices is contained
in the congruence subgroup
\begin{equation*}
\begin{split}
   \left\{\gamma\in\Sp_4(\Z):\gamma\equiv
   \begin{pmatrix} 1 & \ast & \ast & \ast \\
   0 & 1 & \ast & \ast \\ 0 & 0 & 1 & 0 \\
   0 & 0 & \ast & 1 \end{pmatrix} \pmod 5\right\}.
\end{split}
\end{equation*}
The numerical computation in~\cite{CYY} suggests that a similar
phenomenon also occurs in other non-hypergeometric cases. This
naturally leads us to the question whether the Picard--Fuchs
differential equations for Calabi--Yau threefolds are related
to Siegel modular forms in some way.

To be specific, recall the classical result that the solution
${}_2F_1(1/2,1/2;1;z)$ of the Picard--Fuchs differential equation
$$
\theta^2y-\frac z4(2\theta+1)^2y=0
$$
for the family
$$
E_z:y^2=x(x-1)(x-z)
$$
of elliptic curves (i.e., of Calabi--Yau onefolds) satisfies
$$
{}_2F_1\left(\frac12,\frac12;1;\frac{\theta_2^4}{\theta_3^4}\right)
=\theta_3^2,
$$
where $\theta_2(\tau)=\sum_{n\in\Z}e^{\pi i\tau(n+1/2)^2}$ and
$\theta_3(\tau)=\sum_{n\in\Z}e^{\pi i\tau n^2}$ are modular forms of
weight $1/2$. In other words, under a suitable setting, $z$~becomes
a modular function and the holomorphic solution of the
differential equation at $z=0$ becomes a modular form of weight~$1$ on
the congruence subgroup $\Gamma(2)$ of~$\SL_2(\Z)$
(see Section~\ref{s1} for a detailed account of this interpretation).
Likewise, the solution ${}_3F_2(1/4,1/2,3/4;1,1;256z)$ of the
Picard--Fuchs differential equation
\begin{equation}
\theta^3y-4z(4\theta+1)(4\theta+2)(4\theta+3)y=0
\label{3F2}
\end{equation}
for the family
$$
K_z:x_1^4+x_2^4+x_3^4+x_4^4-z^{-1/4}x_1x_2x_3x_4=0
$$
of $K3$ surfaces
(i.e., of Calabi--Yau twofolds) can be interpreted as a modular form of
weight~$2$ on~$\Gamma_0^+(2)$ under a suitable setting. Therefore, one
might be tempted to conjecture that the holomorphic solution of~\eqref{eq: quintic}
at $z=0$ can be interpreted as a Siegel modular form.
The main purpose of the present article is to address this
modularity question.

It turns out that there are several ways to give $\Sp_4$-modular
interpretation, although none of which gives a direct link to Siegel
modular forms. The first of them is given in Section~\ref{s2}, in
which we will show that each fourth order Picard--Fuchs differential
equation for a family of Calabi--Yau threefolds can be associated with
a fifth order linear differential equation, whose holomorphic solution
$w(z)$ at $z=0$, under a suitable formulation, transforms like a
Siegel modular form of weight~$1$ under the action of the monodromy
group. This result can be regarded as a generalization of the
classical Schwarz theory of second order linear differential equations
and automorphic functions.

The formulation of the second modular interpretation is originally due
to A.~Klemm et~al~\cite{ABK}, \cite{HKQ}. Using geometric insights, they showed
that the parameter space of a Picard--Fuchs differential equation of a
family of Calabi--Yau threefolds can be embedded into the Siegel upper
half-space
$$
\H_2=\{\Zeta\in M_2(\C):\trans\Zeta=\Zeta,\ \Im\Zeta>0\},
$$
albeit non-holomorphically. Furthermore, they showed that the image of
the embedding transforms under the action of monodromy in the same way
as~$\H_2$ does under the action of $\Sp_4(\R)$. In Section~\ref{s-klemm},
we will extend this geometric formulation to general fourth
order differential equations with symplectic monodromy. We will show
that under some conditions on the solution space of the differential
equations, we can similarly embed non-holomorphically the parameter
space of the differential equation into~$\H_2$. (If the differential
equation is coming from geometry, the conditions are fulfilled in view
of the argument in~\cite{ABK}, \cite{HKQ}.) Moreover, we will show that if we modify the
function $w(z)$ in the previous paragraph by a certain non-holomorphic
factor, then the resulting function transforms like a Siegel modular
form of weight~$1$ under this formulation.

We stress that, in either interpretation, the functions may not be
related to a true Siegel modular form at all. There are several
reasons for this. One is that the monodromy groups may not be of
finite index in $\Sp_4(\Z)$. (A proof of infiniteness for a particular
example is provided by V.~Pa\c sol in the Appendix. As pointed out to
us by D.~van Straten, the infiniteness of index for the group
generated by the matrices in~\eqref{eq: gen_quintic} might be deduced from
the Margulis--Tits theorem \cite[Chap.~I, \S\,6.6]{Se} modulo some
unproved observation in~\cite{COGP}.) Also, for the
first formulation, the domain on which the function $w$ is defined, is
not always in the Siegel upper half-space. In fact, the discussion in
Section~\ref{s-klemm} shows that for Klemm--et~al's embedding in the second
formulation to be inside~$\H_2$, the domain in the first formulation
has to be disjoint from~$\H_2$. It is a complete mystery how our
functions are related to Siegel modular forms (if they are).

In Section~\ref{s3} we consider the converse problems. Recall that in
the classical Schwarz theory, a second order linear differential
equation with nice monodromy gives rise to a modular function and
a modular form of weight~$1$ and, conversely, for each pair of a modular
form of weight~$1$ and a non-constant modular function, there
associates a second order linear differential equation. In Section~\ref{s3}
we develop an analogous theory for fourth order Picard--Fuchs
differential equations with monodromy inside the symplectic group.


The rest of this article is organized as follows. In Section~\ref{s1}
we review the classical theory of second order linear differential
equations and automorphic functions. This will serve the purpose as a
guideline for our development of a corresponding theory for fourth
order Calabi--Yau differential equations. In Sections~\ref{s2} and~\ref{s-klemm}
we describe the modular interpretations of Calabi--Yau equations
mentioned above. In Section~\ref{s3} we state the converse results
and give some examples in details. The proof of the converse results
will be given in Section~\ref{s4}. Finally, in Section~\ref{s5} we
present some arithmetic observations on a concrete example
of a fifth order Picard--Fuchs linear differential equation.

\section{Modular Picard--Fuchs equations}
\label{s1}

The classical theory concerns with second order linear differential
equations having rational coefficients and regular singularities
at $z=0,\infty$ and some other points. Let us fix such an equation
and assume that its projective monodromy group~$\Gamma$~is a discrete
subgroup of $\SL_2(\mathbb R)$ commensurable with
the full modular group $\SL_2(\mathbb Z)$.
Then we may choose two linearly independent solutions
$u_0=u_0(z)$ and $u_1=u_1(z)$ such that
\begin{equation}
\gamma\colon\begin{pmatrix} u_1 \\ u_0 \end{pmatrix}
\mapsto\chi(\gamma)\begin{pmatrix} a & b \\ c & d \end{pmatrix}
\begin{pmatrix} u_1 \\ u_0 \end{pmatrix}
\qquad\text{for all}\quad
\gamma=\begin{pmatrix} a & b \\ c & d \end{pmatrix}
\in\Gamma,
\label{e01}
\end{equation}
where $\chi(\gamma)$ is a root of unity depending on $\gamma$,
and the image (of a suitably cut $\mathbb C$-plane)
under the multivalued map $\tau(z)=u_1(z)/u_0(z)$
fills the upper half-plane $\H=\{\tau\in\mathbb C:\Im\tau>0\}$.

\begin{Remark}
\label{r1}
If the original differential equation has
the both exponents to be~$0$ at the origin,
the usual choice for a pair of solutions is
$u_0=u_0(z)\in1+z\mathbb C[[z]]$
and $u_1=u_1(z)=\frac{\sqrt D}{2\pi i}(u_0(z)\log z+v(z))$
for some $D\in\mathbb Q$ (normally, $D=1$) and
$v(z)\in z\mathbb C[[z]]$. In this case, the images of
the singular point $z=0$ under the map $\tau(z)$
form a set of cusps.
\end{Remark}

The following is an immediate consequence of our settings
and the Schwarz theory.

\begin{theorem}
\label{tA}
The inverse $z=z(\tau)$ of the map $\tau(z)=u_1(z)/u_0(z)$ is a
\rom(meromorphic\rom) $\Gamma$-modular function \rom(of weight~$0$\rom).
The function $u_0$ viewed as a function of the variable~$\tau$ is a
$\Gamma$-`modular' form of weight~$1$.
\end{theorem}

\begin{Remark}
\label{r2}
We quote the word `modular' since $u_0$ is not necessarily
holomorphic in $\H\cup\{\text{cusps}\}$: it may have branching
(of finite degree) at elliptic points.
Nevertheless, a suitable choice of a positive integer~$N$
gives us a holomorphic $\Gamma$-modular form
$u_0^{N}$ of weight~$N$.
\end{Remark}

\begin{proof}
Indeed, the invariance of $z(\tau)$ under the action of~$\Gamma$
follows from the definition, while from \eqref{e01} we see that
$$
\gamma\colon u_0\mapsto\chi(\gamma)(cu_1+du_0)=\chi(\gamma)(c\tau+d)\cdot
u_0 \qquad\text{for all}\quad \gamma\in\Gamma,
$$
where $\chi(\gamma)$ is a root of unity depending on $\gamma$.
\end{proof}

Therefore, any meromorphic $\Gamma$-modular form of weight~$0$
is an algebraic function of $z(\tau)$, while any meromorphic
$\Gamma$-modular form of weight~$k$ may be represented as
$g(z)u_0^k$, where $g$~is an algebraic function.

\begin{corollary}
The wronskian $W(u_0,u_1)=u_0u_1'-u_0'u_1$
is an algebraic function of~$z$.
\end{corollary}

\begin{proof}[Modular proof]
It follows from \eqref{e01} that
\begin{align*}
\gamma\colon W(u_0,u_1)
\mapsto W(cu_1+du_0,au_1+bu_0)
&=(ad-bc)W(u_0,u_1)
\\
&=W(u_0,u_1)
\qquad\text{for all}\quad \gamma\in\Gamma,
\end{align*}
hence $W(u_0,u_1)$ is a $\Gamma$-modular form
of weight~$0$. This, in particular, implies
the required statement.
\end{proof}

\begin{proof}[Analytic proof]
Let $u''+A(z)u'+B(z)u=0$ denote the original
differential equation satisfied by $u_0$ and~$u_1$.
Then
$$
\frac{\d}{\d z}W(u_0,u_1)
=(u_0u_1'-u_0'u_1)'
=u_0u_1''-u_0''u_1
=-A(z)W(u_0,u_1).
$$
Since the differential equation is of Picard--Fuchs type,
the rational function $A(z)$ may be decomposed
into the following sum of partial fractions:
$$
A(z)=\sum_{k=1}^K\frac{A_k}{z-z_k},
\qquad A_k\in\mathbb C,
$$
where $z_1,\dots,z_K$ are the finite singularities of the
original equation. It remains to integrate the
equation $\d W/\d z=-A(z)W$ to make certain of
the algebraicity of~$W$ as a function of~$z$.
\end{proof}

For further use, write
\begin{equation}
W(u_0,u_1)=\frac1{Cg_0(z)},
\label{e02}
\end{equation}
where $g_0$ is an algebraic function
(that can be explicitly evaluated for the given
second order linear differential equation) and $C\ne0$~is
a certain normalization constant (for instance,
in the case considered in Remark~\ref{r1} the constant~$C$
is usually taken to be $2\pi i$).

This can be also written as
\begin{equation}
\frac{\d\tau}{\d z}
=\frac{W(u_0,u_1)}{u_0^2}
=\frac1{Cg_0u_0^2},
\label{e02a}
\end{equation}
thus showing that $\d z/\d\tau$ is a $\Gamma$-`modular'
form of weight~$2$.

It is interesting that, if we similarly start with
a third order Picard--Fuchs differential equation,
it always comes as the (symmetric) square
of a second order equation of the form described above.
The projective monodromy group of the third order equation is a
discrete subgroup of $O_3(\mathbb R)$ commensurable
with $O_3(\mathbb Z)\simeq\SL_2(\mathbb Z)$. This means
that we are in about the same case as above, and
an (analytic) solution of the differential equation
is a `modular' form of weight~2.

The situation changes drastically if the differential
equation for a function $u_0(z)$ has order
$l\ge4$: it may be the $(l-1)$-st power of
a second order linear differential equation (when we are faced
again to the case discussed above) but there are plenty
of other possibilities in the general case. Note that
in the modular case, i.e., in the case of the reduction
to a second order equation, we have some kind of the inverse
statement:

\begin{theorem}[{\rm folklore;
for three different proofs see \cite[Subsec.~5.4]{Za}}]
\label{tB}
Let $\Gamma$ be a discrete subgroup of $\SL_2(\mathbb R)$
commensurable with $\SL_2(\mathbb Z)$, \ $z=z(\tau)$~a
non-constant meromorphic $\Gamma$-modular form of weight~$0$
and $U_0=U_0(\tau)$~a meromorphic $\Gamma$-modular form of weight~$l-1$.
Then $U_0,\tau U_0,\dots,\tau^{l-1}U_0$ viewed as functions
of~$z$ are linearly independent solutions of an $l$-\textup{th} order
linear differential equation whose coefficients are algebraic
functions of~$z$.
\end{theorem}

In what follows we give analogs of Theorems~\ref{tA} and~\ref{tB}
for linear differential equations of order~4 and~5 in the case
when the projective monodromy group~$\Gamma$ is commensurable with a
discrete subgroup of $\Sp_4(\mathbb Z)\simeq O_5(\mathbb Z)$.

\section
{Differential equations with monodromy $\Sp_4$}
\label{s2}

Consider now a fourth order Picard--Fuchs differential equation
with projective monodromy group $\Gamma\subset\Sp_4(\mathbb R)$
commensurable with a discrete subgroup of $\Sp_4(\mathbb Z)$
(of not necessarily finite index). Assume that there is a point of
maximally unipotent monodromy. Gather its fundamental matrix solution
\begin{equation}
\begin{pmatrix}
u_3 & u_3' & u_3'' & u_3''' \\
u_2 & u_2' & u_2'' & u_2''' \\
u_1 & u_1' & u_1'' & u_1''' \\
u_0 & u_0' & u_0'' & u_0''' \\
\end{pmatrix},
\label{eFM}
\end{equation}
where the basis $u_0,u_1,u_2,u_3$ is chosen
in such a way that
\begin{equation}
W(u_0,u_2)+W(u_1,u_3)=0
\label{e10}
\end{equation}
and the monodromy matrices are in $\Gamma$.
Introduce the functions
$$
w_{jl}=CW(u_j,u_l)=C(u_ju_l'-u_j'u_l), \quad w_{jl}=-w_{lj},
\qquad 0\le j,l\le3,
$$
where $C\ne0$ is a certain normalization constant.
Thanks to~\eqref{e10} we have a linear relation
\begin{equation}
w_{02}+w_{13}=0;
\label{e11}
\end{equation}
there is also a quadratic relation
\begin{equation*}
w_{01}w_{23}+w_{02}w_{13}+w_{03}w_{12}=0,
\end{equation*}
which is tautological in terms of the $u_j$s.
The five linearly independent functions
$$
w_{01}, \quad w_{02}=-w_{13}, \quad w_{03}, \quad w_{12}, \quad w_{23}
$$
form a solution to a fifth order linear differential equation
(the so-called {\it antisymmetric square\/})
with the monodromy conjugate to a subgroup commensurable
to a discrete subgroup of $O_5(\mathbb Z)\simeq\Sp_4(\mathbb Z)$.

\begin{Remark}
\label{r3}
Assuming additionally that the local
exponents of the equation at the origin are all zero
we may always choose the above basis such that
$u_0\in1+z\mathbb C[[z]]$
is the analytic solution at the origin, while
$u_1=\frac1{2\pi i}(u_0\log z+v)$ for some $v\in z\mathbb C[[z]]$.
The inverse $z(t)$ of the mapping $t(z)=u_1(z)/u_0(z)$
is known in the theory of Calabi--Yau threefolds as the \emph{mirror map}.
Note that in this special case we take $C=2\pi i$ and the solution
$w_{01}$ of the above fifth order equation becomes analytic at the origin
after multiplication by~$z$ by our choice of $u_0$ and~$u_1$.
\end{Remark}

\begin{Remark} \label{r3'} We remark that relation \eqref{e10} is a
  feature somehow unique to differential equations with a point of
  maximally unipotent monodromy. According to the general differential
  Galois theory (see \cite{PS}), if the differential Galois group (the
  Zariski closure of the monodromy group) of a fourth order linear
  differential equation with regular singularities is $\Sp_4(\C)$,
  then there exists a basis $u_0,u_1,u_2,u_3$ for the solution
  space such that $W(u_0,u_2)+W(u_1,u_3)$ is invariant under the
  monodromy group. However, this function may not be
  identically zero in general. For instance, the Picard--Fuchs
  differential equation attached to the operator
  $$
    \theta^2\Bigl(\theta-\frac13\Bigr)\Bigl(\theta+\frac13\Bigr)
    -z\Bigl(\theta+\frac12\Bigr)^2\Bigl(\theta+\frac56\Bigr)\Bigl(\theta+\frac76\Bigr)
  $$
  of the family $y^2=x(x-1)(x^3-z)$ of hyperelliptic curves has
  $\Sp_4(\C)$ as its differential Galois group, according to the
  criteria of Beukers and Heckman~\cite{BH}. However, the exterior
  square of this differential equation has order~$6$. In other words,
  $W(u_0,u_2)+W(u_1,u_3)$ is not the zero function. Instead, it is
  equal to $1/(z(1-z))$ (with a suitably chosen basis).

  On the other hand, if a fourth order differential equation with
  $\Sp_4(\C)$ as its differential Galois group has a point of
  maximally unipotent monodromy, then the exterior square also has
  maximally unipotent monodromy at the same singularity. Since the
  highest power of $\log z$ in $W(u_i,u_j)$ is at most~$4$, the
  exterior square has order at most $5$. This forces relation
  \eqref{e10} to exist.
\end{Remark}

Monodromy matrices $\gamma\in\Gamma$ act
by left matrix multiplication:
\begin{equation}
\gamma\colon
\begin{pmatrix} u_3 & u_3' \\ u_2 & u_2' \\
u_1 & u_1' \\ u_0 & u_0' \end{pmatrix}
\mapsto\gamma\cdot\begin{pmatrix} u_3 & u_3' \\ u_2 & u_2' \\
u_1 & u_1' \\ u_0 & u_0' \end{pmatrix}.
\label{e12}
\end{equation}
The matrix
\begin{align}
\label{e12a}
\Tau
&=\begin{pmatrix} u_3 & u_3' \\ u_2 & u_2' \end{pmatrix}
\begin{pmatrix} u_1 & u_1' \\ u_0 & u_0' \end{pmatrix}^{-1}
=\begin{pmatrix} u_3 & u_3' \\ u_2 & u_2' \end{pmatrix}
\begin{pmatrix} -u_0' & u_1' \\ u_0 & -u_1 \end{pmatrix}
\frac1{u_0u_1'-u_0'u_1}
\\
&=\begin{pmatrix} w_{03}/w_{01} & w_{31}/w_{01} \\
w_{02}/w_{01} & w_{21}/w_{01} \end{pmatrix}
\nonumber
\end{align}
is symmetric since $w_{31}=-w_{13}=w_{02}$. Denote
\begin{equation}
\tau_1(z)=\frac{w_{03}}{w_{01}}, \quad
\tau_2(z)=\frac{w_{02}}{w_{01}}=\frac{-w_{13}}{w_{01}}, \quad
\tau_3(z)=\frac{-w_{12}}{w_{01}},
\label{e13}
\end{equation}
hence
\begin{equation}
\Tau=\begin{pmatrix} \tau_1 & \tau_2 \\ \tau_2 & \tau_3 \end{pmatrix},
\qquad \det\Tau=\frac{w_{23}}{w_{01}}.
\label{e14}
\end{equation}
Then we have the standard $\Sp_4$-action on this $2\times2$ symmetric matrix:
\begin{align*}
\gamma\colon\Tau
&=\begin{pmatrix} u_3 & u_3' \\ u_2 & u_2' \end{pmatrix}
\begin{pmatrix} u_1 & u_1' \\ u_0 & u_0' \end{pmatrix}^{-1}
\\
&\mapsto\left(A\begin{pmatrix} u_3 & u_3' \\ u_2 & u_2' \end{pmatrix}
+B\begin{pmatrix} u_1 & u_1' \\ u_0 & u_0' \end{pmatrix}\right)
\left(C\begin{pmatrix} u_3 & u_3' \\ u_2 & u_2' \end{pmatrix}
+D\begin{pmatrix} u_1 & u_1' \\ u_0 & u_0' \end{pmatrix}\right)^{-1}
\\
&=(A\Tau+B)(C\Tau+D)^{-1}=\gamma\Tau
\qquad\text{for}\quad
\gamma=\begin{pmatrix} A & B \\ C & D \end{pmatrix}\in\Gamma.
\end{align*}

Note that the entries of $\Tau=\Tau(z)$ in \eqref{e14}
are algebraically independent over $\mathbb C(z)$. This follows from
the structure of algebraic relations over $\mathbb C(z)$ in the fundamental
matrix solution~\eqref{eFM}; the relations are induced by the differential
Galois group of the starting fourth order equation and the latter group
is isomorphic to the Zariski closure in $\GL_4(\mathbb C)$
of the monodromy group~$\Gamma$ (see, e.g., \cite[Section~5.1]{PS}),
hence is in $\Sp_4(\mathbb C)$.

The multivalued function $\tau=\tau_1(z)$ takes values
in a certain domain $H\subset\nobreak\mathbb C$.
Viewing $\Tau$ as a matrix-valued function of~$\tau$,
we will say that a function
$f(\Tau(\tau))\colon\break H\to\nobreak\mathbb C$ 
is a $\Gamma$-\emph{modular form of weight\/}~$k$
if
$$
f(\gamma\Tau)
=\det(C\Tau+D)^k\cdot f(\Tau)
\qquad\text{for all}\quad
\gamma=\begin{pmatrix} A & B \\ C & D \end{pmatrix}\in\Gamma.
$$

\begin{Theorem}
\label{t1}
The inverse $z=z(\Tau(\tau))$ of the map
$\tau=\tau_1(z)$ in~\eqref{e13},~\eqref{e14} is a
$\Gamma$-modular form of weight~$0$.
Furthermore, the function $w_{01}$ viewed as a function
of~$\Tau=\Tau(\tau)$ is a $\Gamma$-modular form of weight~$1$.
\end{Theorem}

\begin{proof}
The invariance of $z$ under the action of~$\Gamma$
is a consequence of its definition.

Since
\begin{equation} \label{equation: proof of theorem 1}
\begin{split}
\gamma=\begin{pmatrix} A & B \\ C & D \end{pmatrix}
\colon\begin{pmatrix} u_1 & u_1' \\ u_0 & u_0' \end{pmatrix}
&\mapsto C\begin{pmatrix} u_3 & u_3' \\ u_2 & u_2' \end{pmatrix}
+D\begin{pmatrix} u_1 & u_1' \\ u_0 & u_0' \end{pmatrix}
\\
&=(C\Tau+D)\cdot\begin{pmatrix} u_1 & u_1' \\ u_0 & u_0' \end{pmatrix}
\end{split}
\end{equation}
(see \eqref{e12}), taking the determinants of
the both sides (and normalizing by $-C$) we obtain
$$
\gamma\colon w_{01}\mapsto\det(C\Tau+D)\cdot w_{01}.
$$
This shows that $w_{01}$ is a $\Gamma$-modular form
of weight~$1$.
\end{proof}

We now proceed with further computations.
Let us compute $z$-deriva\-tives of $\tau_j$ in~\eqref{e13}.
We have
$$
\frac{\d\tau_1}{\d z}
=\frac{w_{03}'w_{01}-w_{03}w_{01}'}{w_{01}^2}
$$
and for the latter numerator,
\begin{align*}
&
\frac{w_{03}'w_{01}-w_{03}w_{01}'}{(2\pi i)^2}
\\ &\quad
=(u_0u_3''-u_0''u_3)(u_0u_1'-u_0'u_1)
-(u_0u_3'-u_0'u_3)(u_0u_1''-u_0''u_1)
\\ &\quad
=u_0W(u_0,u_1,u_3),
\end{align*}
hence
$$
\frac{\d\tau_1}{\d z}
=\frac{u_0W(u_0,u_1,u_3)}{w_{01}^2}.
$$
In the similar way,
\begin{align*}
\frac{\d\tau_2}{\d z}
&=\frac{w_{01}'w_{13}-w_{01}w_{13}'}{w_{01}^2}
=\frac{-u_1W(u_0,u_1,u_3)}{w_{01}^2}
\\
&=\frac{w_{02}'w_{01}-w_{02}w_{01}'}{w_{01}^2}
=\frac{u_0W(u_0,u_1,u_2)}{w_{01}^2},
\\
\frac{\d\tau_3}{\d z}
&=\frac{w_{01}'w_{12}-w_{01}w_{12}'}{w_{01}^2}
=\frac{-u_1W(u_0,u_1,u_2)}{w_{01}^2}.
\end{align*}
Therefore,
\begin{equation}
\label{equation: tau3'=tau2'2}
\frac{\d\tau_1}{\d z}
:\frac{\d\tau_2}{\d z}
:\frac{\d\tau_3}{\d z}
=1:\biggl(-\frac{u_1}{u_0}\biggr)
:\biggl(-\frac{u_1}{u_0}\biggr)^2
=1:(-t):t^2,
\end{equation}
where $t=t(z)=u_1(z)/u_0(z)$ (the mirror map in the case of Remark~\ref{r3}).

\begin{corollary}
The function
\begin{equation}
U=\det\begin{pmatrix} u_0 & u_2 \\ u_0''' & u_2''' \end{pmatrix}
+\det\begin{pmatrix} u_1 & u_3 \\ u_1''' & u_3''' \end{pmatrix}
\label{e15}
\end{equation}
is an algebraic function of~$z$, say,
\begin{equation}
U=-\frac1{Cg_0(z)}
\label{e15a}
\end{equation}
for an algebraic function $g_0$. Moreover, we have the identity
\begin{equation}
W(u_0,u_1,u_3)
=\frac{u_0}{2\pi ig_0(z)}.
\label{e15b}
\end{equation}
\end{corollary}

\begin{proof}
Differentiating twice equality~\eqref{e10} we obtain
\begin{equation}
\det\begin{pmatrix} u_0 & u_2 \\ u_0'' & u_2'' \end{pmatrix}
+\det\begin{pmatrix} u_1 & u_3 \\ u_1'' & u_3'' \end{pmatrix}
=0
\label{e16}
\end{equation}
and
\begin{equation}
\det\begin{pmatrix} u_0 & u_2 \\ u_0''' & u_2''' \end{pmatrix}
+\det\begin{pmatrix} u_1 & u_3 \\ u_1''' & u_3''' \end{pmatrix}
+\det\begin{pmatrix} u_0' & u_2' \\ u_0'' & u_2'' \end{pmatrix}
+\det\begin{pmatrix} u_1' & u_3' \\ u_1'' & u_3'' \end{pmatrix}
=0.
\label{e17}
\end{equation}
In particular, the last equality gives us the expression
\begin{equation}
U=-\det\begin{pmatrix} u_0' & u_2' \\ u_0'' & u_2'' \end{pmatrix}
-\det\begin{pmatrix} u_1' & u_3' \\ u_1'' & u_3'' \end{pmatrix}.
\label{e18}
\end{equation}
Summing \eqref{e15}, \eqref{e18} and differentiating
we arrive at the equality
\begin{equation}
2U'=\det\begin{pmatrix} u_0 & u_2 \\ u_0'''' & u_2'''' \end{pmatrix}
+\det\begin{pmatrix} u_1 & u_3 \\ u_1'''' & u_3'''' \end{pmatrix}.
\label{e19}
\end{equation}
Using the original fourth order linear differential equation
$u^{(4)}+A(z)u^{(3)}+\dotsb\allowbreak=\nobreak0$, equalities
\eqref{e10}, \eqref{e16} and definition~\eqref{e15}
we finally find from~\eqref{e19} that
\begin{equation}
2U'=-A(z)U.
\label{e20}
\end{equation}
It remains to repeat the argument given in the analytic proof
of Corollary of Theorem~\ref{tA} to get the algebraicity of~$U$
as a function of~$z$.

Expanding the determinant
$$
W(u_0,u_1,u_3)
=\det\begin{pmatrix} u_0 & u_1 & u_3 \\
u_0' & u_1' & u_3'  \\ u_0'' & u_1'' & u_3'' \end{pmatrix}
$$
along the first column and using then
equalities \eqref{e10}, \eqref{e16} and \eqref{e18}
we obtain
\begin{align*}
W(u_0,u_1,u_3)
&=u_0\det\begin{pmatrix} u_1' & u_3' \\ u_1'' & u_3'' \end{pmatrix}
-u_0'\det\begin{pmatrix} u_1 & u_3 \\ u_1'' & u_3'' \end{pmatrix}
+u_0''\det\begin{pmatrix} u_1 & u_3 \\ u_1' & u_3' \end{pmatrix}
\\
&=-u_0\left(U+\det\begin{pmatrix} u_0' & u_2' \\ u_0'' & u_2'' \end{pmatrix}\right)
\\ &\qquad
+u_0'\det\begin{pmatrix} u_0 & u_2 \\ u_0'' & u_2'' \end{pmatrix}
+u_0''\det\begin{pmatrix} u_0 & u_2 \\ u_0' & u_2' \end{pmatrix}
\\
&=-u_0U-W(u_0,u_0,u_2)
=-u_0U
=\frac{u_0}{Cg_0(z)}.
\end{align*}
This proves \eqref{e15b}.
\end{proof}

Recalling that $t(z)=u_1(z)/u_0(z)$ we finally get the formulas
\begin{equation}
\frac{\d t}{\d z}=\frac{w_{01}}{Cu_0^2},
\qquad
\frac{\d\tau_1}{\d z}=\frac{u_0^2}{Cg_0w_{01}^2},
\label{e20a}
\end{equation}
hence
$$
\frac{\d t}{\d z}\cdot\frac{\d\tau_1}{\d z}
=\frac1{C^2g_0w_{01}}.
$$
The resulted formulas may be viewed as analogs of formula~\eqref{e02a}.
(Freely speaking, the product
$(\partial z/\partial t)\cdot(\partial z/\partial\tau_1)$
is a $\Gamma$-`modular' form of weight~1.)

\begin{Remark}
\label{r4}
Some time ago, the first author communicated that
the antisymmetric square of the fifth order
differential equation equals the symmetric square of the corresponding
fourth order differential equation. This became later
the subject of Section~2 in~\cite{Al} (see also Theorem~\ref{theorem:
  converse 3} below).
The theorem easily allows one to deduce that
$w_{03}'w_{01}-w_{03}w_{01}'=C^2u_0W(u_0,u_1,u_3)$
differs from $u_0^2$ by an algebraic function of~$z$,
which is equivalent to the above corollary.
\end{Remark}

\section{A non-holomorphic $\Sp_4$-modularity}
\label{s-klemm}

In this section, we discuss another $\Sp_4$-modular
interpretation for fourth order Picard--Fuchs differential
equations associated with Calabi--Yau threefolds. The formulation is
given in \cite[Section~5]{ABK} and \cite[Section~6]{HKQ}.
Throughout the section we will retain all
the notation used in the previous section.

Let $u_0,u_1,u_2,u_3$ be a basis of the solution space of the
Picard--Fuchs differential equation for a family of Calabi--Yau
threefolds such that the monodromy group with respect to the basis is
$\Gamma\subset\Sp_4(\R)$ (the recipe we follow in the previous section).
Consider, as in \cite[Section~5]{ABK}, the function
$$
F=\frac{u_0u_2+u_1u_3}{2u_0^2}
$$
and let
$$
\tau_{ij}=\frac{\partial^2 F}{\partial u_i\partial u_j}
$$
for $i,j=0,1$. Set
$$
  T=\begin{pmatrix}\tau_{11}&\tau_{10}\\\tau_{01}&\tau_{00}
  \end{pmatrix}.
$$
It is not difficult to see that this $T$ in fact coincides
with our $\Tau$ defined in Section~\ref{s2}.
By geometric consideration, it was shown in~\cite{ABK} that the
imaginary part of~$\Tau$ is indefinite and, thus, is not in the Siegel
upper half-space. Instead, Klemm et~al.\ defined
$$
  \Zeta=\overline\Tau+2i\frac{\Im\Tau\,\trans uu\Im\Tau}
  {u\Im\Tau\,\trans u}, \quad u=(u_1\ u_0),
$$
and showed that $\Zeta$ is in the Siegel upper half-space and transforms
as
$$
  \Zeta\mapsto(A\Zeta+B)(C\Zeta+D)^{-1},
  \qquad\begin{pmatrix}A&B\\C&D\end{pmatrix}\in\Gamma,
$$
under the action of the monodromy group $\Gamma$. Therefore, if we
consider $z$ as a function of $\Zeta$, then $z$ behaves like a Siegel
modular function. In this section, we will extend this result to
general $4$th order linear ordinary differential equations.

\begin{Theorem} \label{theorem: klemm} Let all the notation be given
  as in Section~\textup{\ref{s2}}. Set
  $$
    \phi:=(u_1\ u_0)\Im\Tau\begin{pmatrix}u_1\\ u_0\end{pmatrix}
  $$
  and
  $$
    \Zeta:=\overline\Tau+2i\phi^{-1}\Im\Tau
    \begin{pmatrix}u_1\\ u_0\end{pmatrix}(u_1\ u_0)\Im\Tau.
  $$
  Then, the action of a monodromy
  $\gamma=\left(\begin{smallmatrix}A&B\\C&D\end{smallmatrix}\right)
  \in\Gamma$ on $\Zeta$ is given by
  $$
    \gamma\colon\Zeta\mapsto\gamma\Zeta=(A\Zeta+B)(C\Zeta+D)^{-1}.
  $$
  Also, if we consider $z$ and $\phi\overline w_{01}$ as functions of
  $\Zeta$, then they satisfy
  $$
    z(\gamma\Zeta)=z(\Zeta), \qquad
    \phi(\gamma\Zeta)\overline{w_{01}(\gamma\Zeta)}
   =\det(C\Zeta+D)\phi(\Zeta)\overline{w_{01}(\Zeta)}.
  $$
  Moreover, if the basis $u_0,u_1,u_2,u_3$ satisfies
  \begin{enumerate}
  \item[(i)] $\det(\Im\Tau)<0$ \textup(i.e., $\Im\Tau$ is indefinite\textup) and
  \item[(ii)] $\Im(u_3\overline u_1+u_2\overline u_0)>0$,
  \end{enumerate}
  then $\Zeta$ is contained in the Siegel upper half-space.
\end{Theorem}

\begin{Remark} Note that properties (i) and (ii) assumed in the
  theorem are invariant under the action of the monodromy group
  $\Gamma\subset\Sp_4(\R)$. The first one is very well-known. For the
  second property, we note that the inequality can be written as
  $$
    \frac1{2i}\overline u J\,^tu>0,
  $$
  where $u=(u_3\ u_2\ u_1\ u_0)$ and
  $J=\left(\begin{smallmatrix} O & -E \\ E & O \end{smallmatrix}\right)$,
  \,$O$ and~$E$ stand for the $2\times2$ zero and identity matrices, respectively.
  Now if $\gamma\in\Gamma$, then the action of $\gamma$ gives
  $$
    \frac1{2i}\overline u\,^t\gamma J\gamma\,^tu
   =\frac1{2i}\overline uJ\,^tu>0
  $$
  since $^t\gamma J\gamma=J$.
\end{Remark}

\begin{Lemma} Let all the notation be given as in Theorem
  \textup{\ref{theorem: klemm}}. Then, under the action of
 $\gamma=\left(\begin{smallmatrix}A&B\\C&D\end{smallmatrix}\right)
  \in\Gamma$, we have
  $$
    \phi(\gamma\Zeta)=(u_1\ u_0)\Im\Tau(C\overline\Tau+D)^{-1}(C\Tau+D)
    \begin{pmatrix}u_1\\u_0\end{pmatrix}.
  $$
  Moreover, we have the identity
  $$
    \frac{\phi(\gamma\Zeta)}{\phi(\Zeta)}
   =\frac{\det(C\Zeta+D)}{\det(C\overline\Tau+D)}.
  $$
\end{Lemma}

\begin{proof} From \eqref{equation: proof of theorem 1}, we have
$$
  \gamma:\begin{pmatrix}u_1\\u_0\end{pmatrix}
  \mapsto(C\Tau+D)\begin{pmatrix}u_1\\u_0\end{pmatrix}.
$$
Also, a fundamental property of $\Sp_4(\R)$ states that
$$
  \Im(A\Tau+B)(C\Tau+D)^{-1}=\,\trans{(C\Tau+D)^{-1}}\Im\Tau
  (C\overline\Tau+D)^{-1}.
$$
{}From these two properties, we see that
$$
  \gamma\colon\phi(\Zeta)\mapsto(u_1\ u_0)\,\trans{(C\Tau+D)}\,
  \trans{(C\Tau+D)^{-1}}\Im\Tau(C\overline\Tau+D)^{-1}
  (C\Tau+D)\begin{pmatrix}u_1\\u_0\end{pmatrix}.
$$
This yields the first identity in the lemma. The second identity can
be verified by brute force.
\end{proof}

\begin{proof}[Proof of Theorem~\textup{\ref{theorem: klemm}}]
  The fact that $\Zeta$
  transforms to $(A\Zeta+B)(C\Zeta+D)^{-1}$ can be verified by brute force,
  with the aid of the above lemma. The assertion that
  $z((A\Zeta+B)(C\Zeta+D)^{-1})=z(\Zeta)$ follows from the fact that $z$ is
  invariant under the action of monodromy.  We then combine Theorem
  \ref{t1} with the previous lemma to prove the claim about
  $\phi(\Zeta)\overline{w_{01}(\Zeta)}$. It remains to show that under the two
  conditions in the statements, the function $\Zeta$ is contained in the
  Siegel upper half-space.

  By a direct computation, we find
  $$
    \det(\Im\Zeta)=-\frac{\det(\Im\Tau)}{|\phi(\Zeta)|^2}
    \left((u_1\ u_0)\Im\Tau\begin{pmatrix}\overline u_1\\
    \overline u_0\end{pmatrix}\right)^2.
  $$
  Therefore, if $\Im\Tau$ is indefinite, then $\Im\Zeta$ is either
  positive definite or negative definite, depending on the sign of the
  $(1,1)$-entry of $\Im\Zeta$. Now if we write
  $\Im\Tau=\left(\begin{smallmatrix}a&b\\b&c\end{smallmatrix}\right)$,
  then the $(1,1)$-entry is
  $$
    \left((u_1\ u_0)\Im\Tau\begin{pmatrix}\overline u_1\\
    \overline u_0\end{pmatrix}\right)
    \left((u_1\ u_0)\begin{pmatrix}a^2&ab\\ab&2b^2-ac\end{pmatrix}
    \begin{pmatrix}\overline u_1\\\overline u_0\end{pmatrix}\right).
  $$
  The second factor is a quadratic form of discriminant
  $4a^2(ac-b^2)=4a^2\det(\Im\Tau)$, which by assumption is negative.
  Thus the second factor is positive definite. Therefore, $\Im\Zeta$ is
  positive definite if
  $$
    (u_1\ u_0)\Im\Tau\begin{pmatrix}\overline u_1\\
    \overline u_0\end{pmatrix}>0.
  $$
  Using the relation
  $$
    \Tau\begin{pmatrix}u_1\\u_0\end{pmatrix}
   =\begin{pmatrix}u_3\\u_2\end{pmatrix}
  $$
(cf.~\eqref{e12a}), we find that the condition can be written as
$\Im(u_3\overline u_1+u_2\overline u_0)>0$. This completes the proof of the theorem.
\end{proof}

\section{Converse results}
\label{s3}

In Section \ref{s1}, we have seen that if the projective  monodromy
group of a second order linear differential equation is a discrete
subgroup of $\SL_2(\R)$ commensurable with $\SL_2(\Z)$, then one of
the solutions of differential equations is a modular form of weight
$1$ under a suitable setting. Conversely, Theorem \ref{tB} shows that
if we start out with a modular form $u$ of weight $1$ and a modular
function $z$, then $u$ as a function of~$z$ satisfies a second order
linear differential equation. In this section we develop an
analogous theory in the converse direction for Picard--Fuchs
differential equations of order~$4$.

Given a fourth order Picard--Fuchs differential equation with
symplectic monodromy $\Gamma$ and a point of maximally unipotent
monodromy, in Section \ref{s2} we have seen that if $u_0,u_1,u_2,u_3$
are solutions chosen in a way such that
$$
w_{02}+w_{13}=0,
\qquad
w_{jl}=W(u_j,u_l)=u_j\frac{\d u_l}{\d z}-u_l\frac{\d u_j}{\d z},
$$
then
\begin{enumerate}
\item setting $\tau_1=w_{03}/w_{01}$, $\tau_2=w_{02}/w_{01}$,
$\tau_3=-w_{12}/w_{01}$ and
$\Tau=\left(\begin{smallmatrix}\tau_1&\tau_2\\\tau_2&\tau_3
\end{smallmatrix}\right)$, we have
$$
  \Tau\mapsto (A\Tau+B)(C\Tau+D)^{-1}
$$
under the action of
$\left(\begin{smallmatrix}A&B\\C&D\end{smallmatrix}\right)\in\Gamma$;
\item $w_{01}$ is $\Gamma$-modular form of weight $1$ and $z$ is a
$\Gamma$-modular function in the sense that
\begin{equation*}
z\mapsto z, \qquad
w_{01}\mapsto \det(C\Tau+D)\cdot w_{01}
\end{equation*}
under the action of $\left(\begin{smallmatrix}
A&B\\C&D\end{smallmatrix}\right)\in\Gamma$ (Theorem \ref{t1});
\item $\d\tau_3/\d\tau_1=(\d\tau_2/\d\tau_1)^2$
  (identity~\eqref{equation: tau3'=tau2'2}).
\end{enumerate}
Here we show that if $w(\Tau)$ is a $\Gamma$-modular form of
weight~$1$ and $z(\Tau)$~is a $\Gamma$-modular function in the sense
of~(1),~(2) with $\Tau=\left(\begin{smallmatrix}
\tau_1&\tau_2\\\tau_2&\tau_3\end{smallmatrix}\right)$ satisfying
property~(3) above, then $w$ as a function of~$z$
satisfies a fifth order linear differential equation
(Theorems~\ref{theorem: converse 1} and~\ref{theorem: converse 2}).
Furthermore, there associates a fourth order linear differential
equation whose projective monodromy group contains~$\Gamma$
(Theorems \ref{theorem: converse 3} and~\ref{theorem: converse 4}).

In order for our results to make sense, we shall impose the following
assumptions.

\smallskip
\noindent{\bf Assumptions.} Let $\Gamma$ be a discrete
subgroup of $\Sp_4(\R)$ that may not be commensurable with
$\Sp_4(\Z)$. Let $H$~be a connected domain in the upper
half-plane~$\H$. Assume that $t\colon H\to\C$ is a meromorphic function, and
set, for $\tau\in H$,
$$
  \tau_1=\tau, \qquad
  \tau_2=-\int_{c_1}^{\tau_1}t(\tau)\,\d\tau, \quad
  \tau_3=\int_{c_2}^{\tau_1}t(\tau)^2\,\d\tau,
$$
and
$$
  \Tau(\tau)=\begin{pmatrix}\tau_1&\tau_2\\\tau_2&\tau_3\end{pmatrix}.
$$
Of course, $\tau_2$ and $\tau_3$ are multi-valued, depending on the
choice of the paths of integration. Thus, we assume that for each
residue $r$ of $t(\tau)$, the matrix
$$
  \begin{pmatrix}1&0&0&2\pi ir \\0&1&2\pi ir&0\\
  0&0&1&0\\0&0&0&1\end{pmatrix}
$$
is contained in $\Gamma$, and so is
$$
  \begin{pmatrix}1&0&0&0\\0&1&0&2\pi is\\0&0&1&0\\0&0&0&1\end{pmatrix}
$$
for each residue $s$ of $t(\tau)^2$. Then let $\CC\subset M_2(\C)$ be
the curve defined by
$$
  \CC=\{\Tau(\tau):\tau\in H\}
$$
(with all possible branches of $\Tau$ included), and assume that the
map $\Gamma\times\CC\to M_2(\C)$ given by
\begin{equation} \label{temp1}
  \begin{pmatrix}A&B\\C&D\end{pmatrix}\cdot\Tau(\tau)=
  (A\Tau(\tau)+B)(C\Tau(\tau)+D)^{-1}
\end{equation}
defines a group action of $\Gamma$ on $\CC$.
\medskip

\begin{Theorem}
\label{theorem: converse 1} Under the above assumptions, if
$z,w\colon\CC\to\C$ are non-constant meromorphic functions satisfying
\begin{equation}
\label{temp2}
\begin{aligned}
z(\gamma\Tau)&=z(\Tau),
\\
w(\gamma\Tau)&=\chi(\gamma)\det(C\Tau+D)w(\Tau)
\end{aligned}
\qquad\text{for all}\quad \gamma=\begin{pmatrix} A & B \\ C & D
\end{pmatrix}\in\Gamma,
\end{equation}
where $\chi$ is a character of~$\Gamma$, then $w$, $\tau_1w$, $\tau_2w$,
$\tau_3w$, and $(\tau_1\tau_3-\tau_2^2)w$ viewed as functions of~$z$ are
solutions of a fifth order linear differential equation whose
coefficients are invariant under the substitution of~$\Tau$ by
$\gamma\Tau$ for all $\gamma\in\Gamma$.
\end{Theorem}

\begin{Remark}
\label{r5}
Since the coefficients of Picard--Fuchs differential
equations for Calabi--Yau threefolds are all rational functions, it is
natural to conjecture that the quotient space $\Gamma\backslash\CC$
in Theorem~\ref{theorem: converse 1} can be compactified into a compact
Riemann surface such that \eqref{temp1} and \eqref{temp2} continue
to hold on the compactified curve. We leave problems of this
kind to future study.
\end{Remark}

The fifth order differential equation in Theorem~\ref{theorem: converse 1}
can be made more explicit. This is done in the next theorem.
\medskip

\noindent{\bf Notation.}
Throughout this section, $\theta$ denotes
$z\frac{\d}{\d z}$ while the prime~$'$ stands for the differentiation
with respect to~$\tau$.

\begin{Theorem}
\label{theorem: converse 2}
Let $\Gamma$, $t(\tau)$, $\Tau(\tau)$, $z(\Tau)$, and $w(\Tau)$
be given as in Theorem~\rom{\ref{theorem: converse 1}}. Define
$$
v=\frac{\d t(\tau)}{\d\tau}, \quad
G_1=\frac{\d z/\d\tau}{z}, \quad
G_2=\frac{\d w/\d\tau}{w}, \quad
G_3=\frac{\d v/\d\tau}{v}.
$$
Then the fifth order differential equation in
Theorem~\rom{\ref{theorem: converse 1}} is
\begin{equation*}
\begin{split}
&
\theta^5w+10p_1\theta^4w+(10\theta p_1+35p_1^2+5p_2)\theta^3w
\\ &\quad
+\Bigl(5\theta^2p_1+\frac{15}2\theta p_2+45p_1\theta p_1
+50p_1^3+30p_1p_2\Bigr)\theta^2w
\\ &\quad
+\Bigl(46p_1^2\theta p_1+14p_2\theta p_1+24p_1^4
+2p_3+4p_2^2+11p_1\theta^2p_1
\\ &\quad\qquad
+\frac92\theta^2p_2+\theta^3p_1
+7(\theta p_1)^2+52p_1^2p_2+30p_1\theta p_2\Bigr)\theta w
\\ &\quad
+\bigl(4p_2\theta p_2+9p_1\theta^2 p_2
+7(\theta p_1)(\theta p_2)+26p_1^2\theta p_2+2p_2\theta^2p_1
\\ &\quad\qquad
+20p_1p_2\theta p_1+\theta p_3+\theta^3p_2
+4p_1p_3+24p_1^3p_2+8p_1p_2^2\bigr)w
=0,
\end{split}
\end{equation*}
where
$$
p_1=\frac{2G_1'-G_1(G_2+G_3)}{2G_1^2}, \quad
p_2=\frac{24G_3'-20(G_2+G_3)'+5(G_2+G_3)^2-8G_3^2}{20G_1^2},
$$
and
$$
p_3=\frac{-10G_3'''+40G_3G_3''+21(G_3')^2-54G_3^2G_3'+9G_3^4}{50G_1^4}
$$
are functions invariant under the action of~$\Gamma$.
\end{Theorem}

\begin{Remark}
\label{r7}
Let $Lu=0$ be a Picard--Fuchs differential equation of the family of
  Calabi--Yau threefolds with symplectic monodromy. In practice, $L$~
  has a singular point of maximally unipotent monodromy, which is
  usually assumed to be at $z=0$. Then the functions $u_j$ set up at
  the beginning of this section can be chosen in a way such that $u_0$
  is holomorphic at $z=0$, $u_1=c_1u_0\log z+\cdots$, $u_2=c_2u_0(\log
  z)^3+\cdots$, and $u_3=c_3u_0(\log z)^2+\cdots$. Then the Yukawa
  coupling $K$ for the Calabi--Yau threefolds satisfies
  $$
    K=C\frac{\d^2(u_3/u_0)}{\d(u_1/u_0)^2}
     =C\frac{u_0^3W(u_0,u_1,u_3)}{w_{01}^3}
  $$
  for some constant $C$. In terms of $\tau_j$ given in \eqref{e13},
  this can be expressed as
  $$
  \frac1K=-\frac1C\,\frac{\d^2(w_{13}/w_{01})}{\d(w_{03}/w_{01})^2}
  =-\frac1C\,\frac{\d^2\tau_2}{\d\tau_1^2}.
  $$
  Thus, when the fifth order differential equation in Theorem~\ref{theorem: converse 2}
  arises from the antisymmetric square of~$L$,
  the function~$v$ is actually equal to~$-C/K$. Likewise, we find
  $$
  \frac{\d}{\d\tau_1}=\frac{\d}{\d(w_{03}/w_{01})}
  =\frac CK\,\frac{\d}{\d(u_1/u_0)}=-\frac CK\,\frac{\d}{\d t}.
  $$
  In particular, we have
  $$
  G_1=\frac C{Kz}\,\frac{\d z}{\d(u_1/u_0)}
  =-\frac CK\,\frac{\d z/\d t}z.
  $$
  Now let $f^{(k)}$ denote $\d^k f/\d(u_1/u_0)^k=(-1)^k\d^kf/\d t$. Then we have
  \begin{gather*}
  G_1=\frac{Cz^{(1)}}{Kz}, \qquad G_3=-\frac{CK^{(1)}}{K^2}, \qquad
  G_3'=C^2\frac{2(K^{(1)})^2-K^{(2)}K}{K^4},
  \\
  G_3''=C^3\frac{7K^{(1)}K^{(2)}K-8(K^{(1)})^3-K^2K^{(3)}}{K^6},
  \end{gather*}
  and
  $$
  G_3'''=C^4\frac{11K^{(3)}K^{(1)}K^2+7K^2K^{(2)}-59K^{(2)}(K^{(1)})^2K
  +48(K^{(1)})^4-K^{(4)}K^3}{K^8}.
  $$
  Substituting these expressions into~$p_3$ in Theorem~\ref{theorem: converse 2}
  we find that the function
  $$
  \frac{175(K^{(1)})^4-280K^{(2)}(K^{(1)})^2K+49(K^{(2)}K)^2
  +70K^{(3)}K^{(1)}K^2-10K^{(4)}K^3}{K^4(z^{(1)})^4}
  $$
  should be a function invariant under the action of monodromy.
  Indeed, it was noted in~\cite{LY1} and proved in~\cite{LY2} that the
  function above can be expressed in terms of the coefficients of the
  Picard--Fuchs differential equation.
\end{Remark}

Finally, the last piece of the converse theory shows that under the
assumptions of Theorem~\ref{theorem: converse 1}, there does exist a
fourth order linear differential equation whose projective monodromy
group contains~$\Gamma$.

\begin{Theorem}
\label{theorem: converse 3}
Let all the notations be given as in
Theorems~\rom{\ref{theorem: converse 1}} and~\rom{\ref{theorem: converse 2}}.
Write $w_0=w$, $w_j=\tau_jw$ for $j=1,2,3$, and $w_4=(\tau_1\tau_3-\tau_2^2)w$.
Then the four functions
$$
  \begin{vmatrix}w_0& \theta w_0\\ w_1&\theta w_1\end{vmatrix}^{1/2},
  \quad
  \begin{vmatrix}w_0& \theta w_0\\ w_3&\theta w_3\end{vmatrix}^{1/2},
  \quad
  \begin{vmatrix}w_1& \theta w_1\\ w_4&\theta w_4\end{vmatrix}^{1/2},
  \quad
  \begin{vmatrix}w_3& \theta w_3\\ w_4&\theta w_4\end{vmatrix}^{1/2}
$$
viewed as functions of $z=z(\Tau)$ satisfy a fourth order linear
differential equation whose coefficients are polynomials of $p_i$ and
their derivatives. Moreover, its projective monodromy group contains~$\Gamma$.
\end{Theorem}

\begin{Remark}
\label{r8}
If we let $'$ denote the differentiation with respect
to $\tau_1=\tau$, the four functions in Theorem~\ref{theorem: converse 3}
can be alternatively expressed as
$$
u=\begin{vmatrix}w_0& \theta w_0\\ w_1&\theta w_1\end{vmatrix}^{1/2}, \quad
\tau_2'u, \quad (\tau_1\tau_2'-\tau_2)u, \quad (\tau_2\tau_2'-\tau_3)u,
$$
respectively, up to a sign. To see why this is so, we observe that
$$
\frac{\begin{vmatrix}w_0& \theta w_0\\ w_3&\theta w_3\end{vmatrix}}
{\begin{vmatrix}w_0& \theta w_0\\ w_1&\theta w_1\end{vmatrix}}
=\frac{\d(w_3/w_0)}{\d(w_1/w_0)}=\tau_3'=(\tau_2')^2.
$$
This shows that
$$
\begin{vmatrix}w_0& \theta w_0\\ w_3&\theta w_3\end{vmatrix}^{1/2}
=\pm\tau_2'\begin{vmatrix}w_0& \theta w_0\\ w_1&\theta w_1\end{vmatrix}^{1/2}.
$$
The alternative expressions of
the other two functions can be computed in the same way.
\end{Remark}

\begin{Remark} Note that one can never get the conclusion that the
  projective monodromy group equals to $\Gamma$ in Theorem
  \ref{theorem: converse 3} because the functions $z$ and $w$ may
  actually be modular on a group $\Gamma'$ strictly larger than what
  we assume. (Then the projective monodromy group will also contain
  $\Gamma'$.)
\end{Remark}

Fourth order pullbacks of the fifth order differential equation in
Theorem~\ref{theorem: converse 1} are by no means unique. In practice,
we find that in most cases the following choice has simpler
coefficients in the pullback differential equations.

\begin{Theorem}
\label{theorem: converse 4}
Let all the notations be given as above, and let
$$
g=\exp\left\{-2\int^zp_1\frac{\d z}z\right\}.
$$
Then the four functions
$$
g\begin{vmatrix}w_0& \theta w_0\\ w_1&\theta w_1\end{vmatrix}^{1/2},
\quad
g\begin{vmatrix}w_0& \theta w_0\\ w_3&\theta w_3\end{vmatrix}^{1/2},
\quad
g\begin{vmatrix}w_1& \theta w_1\\ w_4&\theta w_4\end{vmatrix}^{1/2},
\quad
g\begin{vmatrix}w_3& \theta w_3\\ w_4&\theta w_4\end{vmatrix}^{1/2}
$$
satisfy the fourth order linear differential equation
\begin{equation*}
\begin{split}
&\theta^4u+16p_1\theta^3u
+\frac12(187p_1^2+5p_2+38\theta p_1)\theta^2u
\\ &\quad
+\frac12\bigl(22\theta^2p_1+5\theta p_2+294p_1\theta p_1
+472p_1^3+40p_1p_2\bigr)\theta u
\\ &\quad
+\frac1{16}\bigl(-8p_3+9p_2^2+12\theta^2p_2
+40\theta^3p_1+160p_1\theta p_2+124p_2\theta p_1
\\ &\quad\qquad
+4420p_1^2\theta p_1+680p_1\theta^2p_1
+460(\theta p_1)^2+622p_1^2p_2+3465p_1^4\bigr)u
=0.
\end{split}
\end{equation*}
\end{Theorem}

\begin{Remark}
\label{r9}
The reasons why a fourth order pullback exists in the
form given in Theorem~\ref{theorem: converse 3} and why
the particular choice of pullbacks in Theorem~\ref{theorem: converse 4}
tends to have simpler coefficients can be explained as follows.

  In general, given any four differentiable functions $u_0,u_1,u_2,u_3$
  the wronskians have the relation
  $$
    W(W(u_0,u_1),W(u_2,u_3))=-u_0W(u_1,u_2,u_3)+u_1W(u_0,u_2,u_3).
  $$
  This means that the antisymmetric square of the antisymmetric square
  of a linear differential equation is just the tensor product of the
  differential equation with its exterior cube (the differential
  equation satisfied by the wronskians $W(u_1,u_2,u_3)$ for any
  solutions $u_j$ of the original differential equation). Now for a
  fourth order Picard--Fuchs differential equation $Lu=0$ with
  symplectic monodromy, \eqref{e15b} shows that the exterior cube is
  essentially the same as $L$, except for an algebraic factor.
  Therefore, the antisymmetric square of the antisymmetric square of~$L$ is,
  up to an algebraic factor, the symmetric square of $L$. This explains
  the origin of the fourth order pullbacks. The reason why the
  pullback in Theorem~\ref{theorem: converse 4} often has simpler
  coefficients is because it gets rid of the extra algebraic factor
  appearing in the exterior cube of~$L$. See~\cite{Al} for a more
  detailed computation and discussion.
\end{Remark}

The proof of these theorems will be given in the next section. Here we
present some examples first.

\begin{example}
\label{ex1}
Recall that the modular group $\SL_2(\R)$ can
be naturally embedded in $\Sp_4(\R)$ by
$$
\begin{pmatrix} a & b\\ c & d \end{pmatrix}
\mapsto\begin{pmatrix} a & 0 & b & 0 \\ 0 & 1 & 0 & 0 \\
c & 0 & d & 0 \\ 0 & 0 & 0 & 1 \end{pmatrix}.
$$
A less obvious embedding is given by
$$
\iota\colon\begin{pmatrix} a & b \\ c & d \end{pmatrix}
\mapsto\begin{pmatrix}
a^2d+2abc & -3a^2c & abd+\frac12b^2c & \frac12b^2d \\
-a^2b & a^3 & -\frac12ab^2 & -\frac16b^3 \\
4acd+2bc^2 & -6ac^2 & ad^2+2bcd & bd^2 \\
6c^2d & -6c^3 & 3cd^2 & d^3 \end{pmatrix}.
$$
The origin of this embedding can be explained as follows.

Let $\Gamma$ be a congruence subgroup of $\SL_2(\Z)$
and $f(\tau)$~a modular form of weight~$3$ on~$\Gamma$ with character~$\chi$.
Let $z(\tau)$~be a modular function on~$\Gamma$. Then we have
the equality of wronskians,
$$
W(f,-\tfrac16\tau^3f)
=f^2W(1,-\tfrac16\tau^3)=-f^2\cdot\frac12\tau^2\frac{\d\tau}{\d z}
=-f^2W(\tau,\tfrac12\tau^2)
=-W(\tau f,\tfrac12\tau^2f),
$$
hence if we let $\wh\gamma$ denote the $4\times 4$ matrix satisfying
$$
  \begin{pmatrix} \frac12\tau^2f \\ -\frac16\tau^3f \\ \tau f \\ f
  \end{pmatrix}|\gamma
  =\wh\gamma\cdot\begin{pmatrix} \frac12\tau^2f \\ -\frac16\tau^3f \\
    \tau f \\ f \end{pmatrix} \qquad\text{for}\quad
  \gamma=\begin{pmatrix} a & b \\ c & d  \end{pmatrix}\in\Gamma
$$
then, up to a numerical scalar, $\wh\gamma$~is in the symplectic group.
Indeed, a direct computation then shows that, up to the character
$\chi$, the matrix $\wh\gamma$ is $\iota\gamma$.

Now let $t(\tau)=\tau$. Set
$$
\tau_1=\tau, \quad
\tau_2=-\int_0^{\tau}t(\tau)\,\d\tau=-\frac12\tau^2, \quad
\tau_3=\int_0^{\tau}t(\tau)^2\d\tau=\frac13\tau^3
$$
and
\begin{equation}
\label{equation: example 1 1}
\Tau(\tau)=\begin{pmatrix} \tau_1 & \tau_2 \\ \tau_2 & \tau_3 \end{pmatrix}
=\begin{pmatrix} \tau & -\frac12\tau^2 \\ -\frac12\tau^2 & \frac13\tau^3 \end{pmatrix}.
\end{equation}
(This is exactly the choice imposed by formulas~\eqref{e13} and~\eqref{e14}.)
For $\gamma=\bigl(\begin{smallmatrix} a & b \\ c & d \end{smallmatrix}\bigr)\in\Gamma$,
if we write $\iota\gamma=\bigl(\begin{smallmatrix} A & B \\ C & D \end{smallmatrix}\bigr)$
then we have
$$
(A\Tau+B)(C\Tau+D)^{-1}
=\begin{pmatrix} \gamma\tau & -\frac12(\gamma\tau)^2 \\
-\frac12(\gamma\tau)^2 & \frac13(\gamma\tau)^3 \end{pmatrix}
=\Tau(\gamma\tau),
$$
where $\gamma\tau=(a\tau+b)/(c\tau+d)$. Thus, the mapping
$$
(\iota\gamma,\Tau)\mapsto(A\Tau+B)(C\Tau+D)^{-1}
$$
defines a group action of $\iota(\Gamma)$ on the set
$\CC=\{\Tau(\tau):\tau\in\H\}$.

Moreover, we have
$$
\det(C\Tau+D)=(c\tau+d)^4(ad-bc)=(c\tau+d)^4.
$$
Thus, if $w(\tau)$ is a modular form of weight~$4$ on~$\Gamma$, then
Theorem~\ref{theorem: converse 1} implies that $w$, $\tau_1w=\tau w$,
$\tau_2w=-\frac12\tau^2w$, $\tau_3w=\frac13\tau^3w$, and
$(\tau_1\tau_3-\tau_2^2)w=\frac1{12}\tau^4w$, as functions of~$z$, satisfy a
fifth order linear differential equation with algebraic functions as
coefficients, in accordance with Theorem~\ref{tB}.
\end{example}

\begin{Remark}
\label{r10}
Note that in the above example we have
$\det\Im\Tau(\tau)=-(\Im\tau)^4/3$. Thus, the
curve $\CC=\{\Tau(\tau):\tau\in\H\}$ is not contained in
the Siegel upper half-space. In other words, the functions $z$ and $w$
in Theorem~\ref{theorem: converse 1} may not be related to
Siegel modular functions and modular forms at all.
\end{Remark}

\begin{example}
\label{ex2}
Consider the Picard--Fuchs differential equation~\eqref{eq: quintic}
for the quintic threefolds. Let
\begin{gather*}
y_0=1+120z+113400z^2+\cdots,
\qquad y_1=\frac1{2\pi i}(y_0\log z+g_1),
\\
y_2=\frac1{(2\pi i)^2}\biggl(y_0\frac{\log^2z}2+g_1\log z+g_2\biggr),
\\
y_3=\frac1{(2\pi i)^3}\biggl(y_0\frac{\log^3z}6+g_1\frac{\log^2z}2+g_2\log z+g_3\biggr)
\end{gather*}
be the (normalized) Frobenius basis at $z=0$. In~\cite{CYY} we showed that with
respect to the ordered basis
\begin{gather*}
u_3=5y_2+\frac52y_1-\frac{25}{12}y_0,
\quad
u_2=-5y_3-\frac{25}{12}y_1+\frac{200\zeta(3)}{(2\pi i)^3}y_0,
\quad
u_1=y_1, \quad u_0=y_0,
\end{gather*}
the monodromy matrices around $z=0$ and $z=1/3125$ are
$$
  \begin{pmatrix}1&0&5&5\\-1&1&0&-5\\0&0&1&1\\0&0&0&1\end{pmatrix}
  \quad\text{and}\quad
  \begin{pmatrix}1&0&0&0\\0&1&0&0\\0&0&1&0\\0&1&0&1\end{pmatrix},
$$
respectively. (Note that this is obtained from the collection~\eqref{eq: gen_quintic}
given in the introduction by reordering the basis.) Write
$w_{jl}=W(u_j,u_l)$, $\tau_1=w_{03}/w_{01}$, and
$\tau_2=w_{02}/w_{01}$. We find
\begin{gather*}
w_{01}=\frac1z+1010+1861650z+4119140000z^2+9959217231250z^3+\cdots,
\\
\tau_1=\frac1{2\pi i}\left(5\log z+5\pi i+6725z+\frac{16482625}2z^2
+\frac{44704818125}3z^3+\cdots\right),
\end{gather*}
and
$$
\tau_2=-\frac{\tau_1^2}{10}+\frac{\tau_1}2+\frac1{(2\pi i)^2}
\left(\frac{65}6\pi^2+2875z+\frac{17038125}4z^2+\frac{151564765625}{18}z^3+\cdots\right).
$$
Then the functions $v$ and $G_j$ in Theorem~\ref{theorem: converse 2}
have the $z$-expansions
$$
v=-\frac15+115z+217500z^2+471493250z^3+1103069708750z^4+\cdots
$$
(notice that the Yukawa coupling has the $z$-expansion
$5+2875z+7090625z^2+18991003125z^3+\cdots$, which is exactly $-1/v$),
\begin{align*}
G_1&=2\pi i\left(\frac15-269z-297500z^2-501290000z^3-1001288510000z^4-\cdots\right),
\\
G_2&=2\pi i\left(-\frac15+471z+566450z^2+1023038500z^3+2170808632500z^4+\cdots\right),
\\
G_3&=2\pi i(-115z-346450z^2-982613500z^3-2787375077500z^4-\cdots).
\end{align*}
We find
\begin{gather*}
p_1=\frac{1-6250z}{2(1-3125z)}, \qquad
p_2=\frac{1-17000z+37500000z^2}{(1-3125z)^2},
\\
p_3=\frac{5z(46+509375z+156250000z^2)}{(1-3125z)^4},
\end{gather*}
and the fifth order differential equation in Theorem~\ref{theorem: converse 2}
for the functions $\wt w=zw$ is
\begin{equation*}
\begin{split}
&\theta^5\wt w-5z(2\theta+1)(625\theta^4+1250\theta^3+1500\theta^2+875\theta+202)\wt w
\\ &\qquad
+5^5z^2(5\theta+3)(5\theta+4)(5\theta+5)(5\theta+6)(5\theta+7)\wt w=0.
\end{split}
\end{equation*}
(We normalize the functions $w\mapsto zw$ in order to have the local exponents
zero with multiplicity~5 at $z=0$.)
Using the {\tt Maple} command {\tt exterior\_power}, we find that this
is indeed the antisymmetric square
of the differential equation~\eqref{eq: quintic}.
\end{example}

\begin{example}
\label{ex3}
Consider the differential equation
\begin{equation}
\theta^5y-32z(2\theta+1)^5y=0
\label{e23}
\end{equation}
with singularities at the points $z=0$, $1/1024$, and~$\infty$.
Let $y_0,y_1,y_2,y_3,y_4$ denote its (normalized) Frobenius basis at $z=0$,
\begin{gather*}
y_0=f_0(z),
\qquad
y_1=\frac1{2\pi i}\bigl(f_0(z)\log z+f_1(z)\bigr),
\\
y_2=\frac1{(2\pi i)^2}\biggl(f_0(z)\frac{\log^2z}2+f_1(z)\log z+f_2(z)\biggr),
\\
y_3=\frac1{(2\pi i)^3}\biggl(f_0(z)\frac{\log^3z}{3!}+f_1(z)\frac{\log^2z}2
+f_2(z)\log z+f_3(z)\biggr),
\\
y_4=\frac1{(2\pi i)^4}\biggl(f_0(z)\frac{\log^4z}{4!}
+f_1(z)\frac{\log^3z}{3!}+f_2(z)\frac{\log^2z}2+f_3(z)\log z+f_4(z)\biggr).
\end{gather*}
Following Beukers' argument \cite[Sections~3 and~4]{AZ},
we can show that these functions satisfy
$$
y_0y_4-y_1y_3+\frac12y_2^2=0, \qquad
(\theta y_0)(\theta y_4)-(\theta y_1)(\theta y_3)+\frac12(\theta y_2)^2=0
$$
Thus, up to conjugation, the monodromy group is contained in the
orthogonal group~$O_5$. Using the method in~\cite{CYY} we can prove
that, relative to the ordered basis $y_4,y_3,y_2,y_1,y_0$,
the monodromy matrices around $z=0$ and $z=1/1024$ are
$$
  \MF{11{1/2}{1/6}{1/24}}{011{1/2}{1/6}}{0011{1/2}}
  {00011}{00001}, \quad
  \MF{{a^2}0{-ab}{(1-a^2)x}{-b^2/2}}
  {{-c^2x/2}1{-acx}{c^2x^2/2}{-(1-a^2)x}}
  {{-ac}0{1-2a^2}{acx}{-ab}}{00010}
  {{-c^2/2}0{-ac}{c^2x/2}{a^2}},
$$
where $a=5/6$, $b=11/144$, $c=8$, and $x=10\zeta(3)/(2\pi i)^3$
(see \cite[Theorem~3]{CYY}). Set
$$
  \begin{pmatrix}w_4\\ w_3\\ w_2\\ w_1\\ w_0\end{pmatrix}
  =\begin{pmatrix}
   32 & 0 & 20/3 & -32x & -25/36 \\
   0 & 8 & 0 & 0 & -8x \\
   0 & 0 & -4 & 0 & 5/6 \\
   0 & 0 & 0 & 4 & 0 \\
   0 & 0 & 0 & 0 & 1 \end{pmatrix}
  \begin{pmatrix} y_4 \\ y_3 \\ y_2 \\ y_1 \\ y_0 \end{pmatrix}.
$$
With respect to this new basis, the matrices become
$$
  \MF{14{-4}38}{01{-2}13}{001{-1}{-2}}{00014}{00001}
  \quad\text{and}\quad
  \MF{0000{-1}}{01000}{00100}{00010}{{-1}0000}.
$$
Also, $w_j$ satisfy
$$
  w_0w_4-w_1w_3+w_2^2=0, \qquad
 (\theta w_0)(\theta w_4)-(\theta w_1)(\theta w_3)+(\theta w_2)^2=0.
$$
Now for $j=1,\dots,4$, let $\tau_j=w_j/w_0$, and let $\tau_j'$
denote the derivative of $\tau_j$ with respect to $\tau=\tau_1$. From the
above relations we deduce that
$$
\tau_4=\tau_1\tau_3-\tau_2^2, \qquad \tau_3'=(\tau_2')^2.
$$
Set $\Tau(\tau)=\Tau(\tau_1)=\bigl(\begin{smallmatrix}\tau_1&\tau_2\\
\tau_2&\tau_3\end{smallmatrix}\bigr)$.
Around $z=0$, we have $w_0\mapsto w_0$ and
\begin{equation*}
\begin{split}
\Tau(\tau)
&\mapsto \begin{pmatrix} \tau_1+4 & \tau_2-\tau_1-2 \\
\tau_2-\tau_1-2 & \tau_3-2\tau_2+\tau_1+3 \end{pmatrix}
\\
&=\left(\begin{pmatrix} 1 & 0 \\ -1 & 1 \end{pmatrix}\Tau(\tau)
+\begin{pmatrix} 4 & 2 \\ -2 & 1 \end{pmatrix}\right)
\left(\begin{pmatrix} 0 & 0 \\ 0 & 0 \end{pmatrix}\Tau(\tau)
+\begin{pmatrix} 1 & 1 \\ 0 & 1 \end{pmatrix}\right)^{-1}.
\end{split}
\end{equation*}
Around $z=1/1024$, we have
$$
w_0\mapsto-\tau_4w_0=-(\tau_1\tau_3-\tau_2^2)w_0
=-\det\left(\begin{pmatrix} 0 & -1 \\ 1 & 0 \end{pmatrix}\Tau(\tau)
+\begin{pmatrix} 0 & 0 \\ 0 & 0 \end{pmatrix}\right)w_0
$$
and
\begin{equation*}
\begin{split}
\Tau(\tau)
&\mapsto-\frac1{\tau_4}
\begin{pmatrix} \tau_1 & \tau_2\\ \tau_2 & \tau_3 \end{pmatrix}
\\
&=\left(\begin{pmatrix} 0 & 0 \\ 0 & 0 \end{pmatrix}\Tau(\tau)
+\begin{pmatrix} 0 & 1 \\ -1 & 0 \end{pmatrix}\right)
\left(\begin{pmatrix} 0 & -1 \\ 1 & 0 \end{pmatrix}\Tau(\tau)
+\begin{pmatrix} 0 & 0 \\ 0 & 0 \end{pmatrix}\right)^{-1}.
\end{split}
\end{equation*}
Thus, letting $\Gamma$ be the subgroup of $\Sp_4(\Z)$ generated by
\begin{equation}
\gamma_0=\M{1042}{{-1}1{-2}1}{0011}{0001}
\quad\text{and}\quad
\gamma_1=\M{0001}{00{-1}0}{0{-1}00}{1000}=\gamma_1^{-1},
\label{gamma01}
\end{equation}
$w_0(\Tau)$ is a $\Gamma$-modular form of weight~$1$. The functions
$v$ and $G_j$ in Theorem~\ref{theorem: converse 2} have the $z$-expansions
\begin{equation*}
\begin{split}
    v&=1+160z+132320z^2+115614720z^3+104797147360z^4+\cdots, \\
  G_1&=2\pi i(1-160z-54880z^2-29946880z^3-19691390560z^4-\cdots), \\
  G_2&=2\pi i(32z+9408z^2+4805632z^3+3045669248z^4+\cdots), \\
  G_3&=2\pi i(160z+213440z^2+240399360z^3+259173946240z^4+\cdots).
\end{split}
\end{equation*}
We find
$$
  p_1=\frac{-256z}{1-1024z}, \quad
  p_2=\frac{65536z^2}{(1-1024z)^2}, \quad
  p_3=\frac{-32z-163840z^2-33554432z^3}{(1-1024z)^4}.
$$
Of course, the fifth order differential equation in
Theorem~\ref{theorem: converse 2} is just the original hypergeometric
differential equation, while the fourth order differential equation
in Theorem~\ref{theorem: converse 4} is
\begin{equation}
\theta^4-16z(128\theta^4+256\theta^3+304\theta^2+176\theta+39)
+2^{20}z^2(\theta+1)^4.
\label{e23pullback}
\end{equation}
\end{example}

Example~\ref{ex3} is our basic example for further illustrations.
In Section~\ref{s5} we discuss some arithmetic observations around
this example, while Theorem~\ref{th:inf} in the Appendix shows that
$\Gamma$~is not of finite index in~$\Sp_4(\mathbb Z)$.

\section{Proof of the converse theorems}
\label{s4}

\noindent{\bf Notation.}
For
$$
\gamma
=\begin{pmatrix} A & B \\ C & D \end{pmatrix}
=\begin{pmatrix} a_{11} & a_{12} & a_{13} & a_{14} \\
  a_{21} & a_{22} & a_{23} & a_{24} \\
  a_{31} & a_{32} & a_{33} & a_{34} \\
  a_{41} & a_{42} & a_{43} & a_{44} \end{pmatrix}
  \in\Gamma\subset\Sp_4(\R)
$$
and $\tau\in H$ we write $\phantom{\fa}$
\begin{equation}
\label{equation: abcd}
C\Tau+D=\begin{pmatrix} \fa & \fb \\ \fc & \fd \end{pmatrix}
=\begin{pmatrix}
a_{31}\tau_1+a_{32}\tau_2+a_{33} &
a_{31}\tau_2+a_{32}\tau_3+a_{34} \\
a_{41}\tau_1+a_{42}\tau_2+a_{43} &
a_{41}\tau_2+a_{42}\tau_3+a_{44}
\end{pmatrix}.
\end{equation}
The notation $\gamma\tau$ will represent the $(1,1)$-entry of
$\gamma\Tau(\tau)=(A\Tau+B)(C\Tau+D)^{-1}$.
More generally, for a function $g(\tau)$ of~$\tau$
we often write $\gamma g(\tau)$ or $g|\gamma$ in place of $g(\gamma\tau)$.

For ease of notation, we let $M^\adj$
denote the adjugate of a square matrix~$M$,
i.e., $M^\adj$ is the square matrix such that
$M^\adj M=(\det M)\operatorname{Id}$.
\medskip

We start out by doing some elementary computation.

\begin{Lemma}
\label{lemma: Z derivative}
Let $\Tau=\Tau(\tau)$ be defined as in Theorem~\rom{\ref{theorem: converse 1}}.
For $\gamma=\bigl(\begin{smallmatrix} A & B \\ C & D \end{smallmatrix}\bigr)\in\Gamma$,
we have
\begin{gather}
\label{equation: Z derivative 1}
\d(\gamma\Tau)
=\frac{\trans{(C\Tau+D)^\adj}\,\d\Tau\,(C\Tau+D)^\adj}{\det(C\Tau+D)^2},
\\
\label{equation: Z derivative 2}
\frac{\d\gamma\tau}{\d\tau}
=\frac{(\fc t+\fd)^2}{\det(C\Tau+D)^2},
\end{gather}
\begin{equation}
\label{equation: Z derivative 3}
\frac{\d\Tau(\tau)}{\d\tau}|\gamma
=\frac{\trans{(C\Tau+D)^\adj}\,(\d\Tau/\d\tau)\,(C\Tau+D)^\adj}{(\fc t+\fd)^2},
\end{equation}
and
\begin{equation} \label{equation: t}
  t(\gamma\tau)=\frac{\fa t(\tau)+\fb}{\fc t(\tau)+\fd}
  \text{ for all }\gamma\in\Gamma\text{ and }\tau\in\H.
\end{equation}
\end{Lemma}

\begin{proof}
Identity~\eqref{equation: Z derivative 1} is simply a restatement of
the basic property
$$
\d(A\Tau+B)(C\Tau+D)^{-1}=\trans{(C\Tau+D)^{-1}}\,\d\Tau\,(C\Tau+D)^{-1}
$$
of the symplectic group.

To prove \eqref{equation: Z derivative 2}, we compare the
$(1,1)$-entries of the two sides of~\eqref{equation: Z derivative 1}.
The $(1,1)$-entry of the left-hand side is $\d\gamma\tau/\d\tau$,
while the numerator of the right-hand side is
\begin{equation}
\label{equation: Z derivative 4}
\begin{split}
\trans{(C\Tau+D)}^\adj\frac{\d\Tau(\tau)}{\d\tau}(C\Tau+D)^\adj
&=\begin{pmatrix} \fd & -\fc \\ -\fb & \fa \end{pmatrix}
  \begin{pmatrix} 1 & -t \\ -t & t^2 \end{pmatrix}
  \begin{pmatrix} \fd & -\fb \\ -\fc & \fa \end{pmatrix} \\
&=\begin{pmatrix} (\fc t+\fd)^2 & -(\fa t+b)(\fc t+\fd) \\
  -(\fa t+b)(\fc t+\fd) & (\fa t+\fb)^2\end{pmatrix},
\end{split}
\end{equation}
whose $(1,1)$-entry is $(\fc t+\fd)^2$. This proves
identity~\eqref{equation: Z derivative 2}.

Finally, \eqref{equation: Z derivative 3} follows from the
first two, and \eqref{equation: t} follows from \eqref{equation: Z
  derivative 3} and \eqref{equation: Z derivative 4}. This completes
the proof.
\end{proof}

\begin{Lemma} \label{lemma: K}
Let $v(\tau)=\d t/\d\tau$. Then
\begin{equation}
\label{equation: K}
v(\gamma\tau)
=\frac{\det(C\Tau+D)^3}{(\fc t+\fd)^4}v(\tau).
\end{equation}
\end{Lemma}

\begin{proof}
Differentiating~\eqref{equation: t}
we obtain
$$
v(\gamma\tau)\frac{\d\gamma\tau}{\d\tau}
=\frac{(\fa't+\fb')(\fc t+\fd)-(\fa t+\fb)(\fc't+\fd')}{(\fc t+\fd)^2}
+\frac{\fa\fd-\fb\fc}{(\fc t+\fd)^2}v(\tau).
$$
We then observe that from~\eqref{equation: abcd} we have
$\fa\fd-\fb\fc=\det(C\Tau+D)$ and $\fa't+\fb'=\fc't+\fd'=0$.
Then from~\eqref{equation: Z derivative 2} in Lemma~\ref{lemma: Z
  derivative} we obtain~\eqref{equation: K}.
\end{proof}

\begin{Lemma}
\label{lemma: invariance of p}
Let $G_j$ and $p_j$, $j=1,2,3$, be defined as in
Theorem~\rom{\ref{theorem: converse 2}}.
Then $p_j(\gamma\tau)=p_j(\tau)$ for all $\gamma\in\Gamma$.
\end{Lemma}

\begin{proof}
For the sake of convenience, for
$\gamma=\bigl(\begin{smallmatrix} A & B \\ C & D
\end{smallmatrix}\bigr)\in\Gamma$
and $\Tau=\Tau(\tau)\in\CC$, set
$$
g(\gamma,\Tau)=\fc t+\fd,
\qquad
h(\gamma,\Tau)=\det(C\Tau+D).
$$
Taking the logarithmic derivatives of the two sides of~\eqref{temp2}
and~\eqref{equation: K} with respect to~$\tau$ and
then using \eqref{equation: Z derivative 2}, we obtain
\begin{gather}
\label{equation: p 1}
G_1(\gamma\tau)=\frac{h^2}{g^2}G_1(\tau),
\\
\label{equation: p 2}
G_2(\gamma\tau)=\frac{h'h}{g^2}+\frac{h^2}{g^2}G_2(\tau),
\end{gather}
and
\begin{equation}
\label{equation: p 3}
G_3(\gamma\tau)
=3\frac{h'h}{g^2}-4\frac{g'h^2}{g^3}+\frac{h^2}{g^2}G_3(\tau).
\end{equation}

Differentiating \eqref{equation: p 1} with respect to~$\tau_1$
again, we have
$$
G_1'(\gamma\tau)\frac{\d\gamma\tau}{\d\tau}
=2\frac{h'h}{g^2}G_1(\tau)-2\frac{g'h^2}{g^3}G_1(\tau)+\frac{h^2}{g^2}G_1'(\tau).
$$
{}From this, \eqref{equation: Z derivative 2}, \eqref{equation: p 2}, and \eqref{equation: p 3}
we deduce that
$$
\biggl(G_1'-\frac{G_1(G_2+G_3)}2\biggr)|\gamma
=\frac{h^4}{g^4}\biggl(G_1'-\frac{G_1(G_2+G_3)}2\biggr).
$$
It follows that $p_1=(G_1'-G_1(G_2+G_3)/2)/G_1^2$ is invariant under
the action of~$\Gamma$.

We next prove that $p_2$ is invariant under~$\Gamma$. By a direct
computation, we find
\begin{equation*}
\begin{split}
&
\bigl(24G_3'-20(G_2+G_3)'+5(G_2+G_3)^2-8G_3^2\bigr)|\gamma
\\ &\qquad
-\frac{h^4}{g^4}\bigl(24G_3'-20(G_2+G_3)'+5(G_2+G_3)^2-8G_3^2\bigr)
\\ &\quad
=8\frac{h^3}{g^5}(2g'h'+2g'hG_3-2g''h-h''g)
\\ &\quad
=8\frac{h^3}{g^5}(a_{41}a_{33}+a_{42}a_{34}-a_{31}a_{43}-a_{32}a_{44})gv.
\end{split}
\end{equation*}
The property $C\,\trans{D}=D\,\trans{C}$ of the symplectic
matrix~$\gamma$ implies
$a_{41}a_{33}+a_{42}a_{34}-a_{31}a_{43}-a_{32}a_{44}=\nobreak0$.
It follows that $p_2$ is invariant under the action of~$\Gamma$.

The invariance of $p_3$ under $\Gamma$ can be proved in the same way.
We repeatedly use $2g'h'+2g'hG_3-2g''h-h''g=0$ just shown
and~\eqref{equation: Z derivative 2}.
The details are too complicated to be presented here.
\end{proof}

\begin{proof}[Proof of Theorems~\rom{\ref{theorem: converse 1}}
  and~\rom{\ref{theorem: converse 2}}]
We first give a conceptual proof of Theorem~\ref{theorem: converse 1},
which may be of interest independent of our later proof of
Theorem~\ref{theorem: converse 2}.

For convenience, we let $\bw$ denote the column vector
$\trans{((\tau_1\tau_3-\tau_2^2)w,\tau_3w,\tau_2w,\tau_1w,w)}$.
Given $\gamma=\bigl(\begin{smallmatrix} A & B \\ C & D \end{smallmatrix}\bigr)\in\Gamma$,
we have
\begin{align*}
\bigl(w(\Tau)\Tau\bigr)|\gamma
&=\chi(\gamma)\det(C\Tau+D)w(\Tau)(A\Tau+B)(C\Tau+D)^{-1}
\\
&=\chi(\gamma)w(\Tau)(A\Tau+B)(C\Tau+D)^\adj
\\
&=\chi(\gamma)w(\Tau)\bigl((\det\Tau)AC^\adj+A\Tau D^\adj
+B\Tau^\adj C^\adj+BD^\adj\bigr).
\end{align*}
This shows that there is a matrix $\wh\gamma$ in~$M_5(\R)$ such that
$$
\bw|\gamma=\chi(\gamma)\wh\gamma\cdot\bw.
$$
Since $z(\Tau)$ is assumed to be invariant under the substitution
$\Tau\mapsto\gamma\Tau$, the same matrix $\wh\gamma$ also satisfies
$$
\theta^j\bw|\gamma
=\chi(\gamma)\wh\gamma\cdot\theta^j\bw
\qquad\text{for}\quad j=1,2,\dots\,.
$$
Therefore, the coefficients $r_j(\Tau)$ in the linear dependence
$$
\theta^5\bw=\sum_{j=0}^4r_j(\Tau)\theta^j\bw
$$
must be invariant under the action of~$\gamma$. This proves
Theorem~\ref{theorem: converse 1}. We now prove
Theorem~\ref{theorem: converse 2}.

Setting $G_2/G_1=\lambda$ we have
\begin{equation}
\label{equation: theta w}
\theta w=t\frac{\d w/\d\tau}{\d t/\d\tau}=w\lambda
\end{equation}
and
\begin{align}
\label{equation: theta lambda}
\theta\lambda
&=\frac1{G_1}\biggl(\frac{G_2}{G_1}\biggr)'
\\
&=\frac{20G_2'-4G_3'-5G_2^2-10G_2G_3+3G_3^2}{20G_1^2}
\nonumber\\ &\qquad
-\frac{G_2(2G_1'-G_1(G_2+G_3))}{2G_1^3}
-\frac{G_2^2}{4G_1^2}+\frac{4G_3'-3G_3^2}{20G_1^2}
\nonumber\\
&=-p_2-p_1\lambda-\frac14\lambda^2+\mu,
\nonumber
\end{align}
where we let
\begin{equation} \label{equation: p1 p2}
p_1=\frac{2G_1'-G_1(G_2+G_3)}{G_1^2}, \quad
p_2=\frac{4G_3'-20G_2'+5G_2^2+10G_2G_3-3G_3^2}{G_1^2},
\end{equation}
and
$$
\mu=\frac{4G_3'-3G_3^2}{20G_1^2}.
$$
Note that $p_1$ and $p_2$ are invariant under the action of~$\Gamma$
by Lemma~\ref{lemma: invariance of p}.

Furthermore, by a similar argument we have
\begin{equation}
\label{equation: theta mu}
\theta\mu
=\frac1{G_1}\biggl(\frac{4G_3'-3G_3^2}{20G_1^2}\biggr)
=-2p_1\mu-\lambda\mu+\nu,
\end{equation}
where
$$
\nu=\frac{4G_3''-10G_3'G_3+3G_3^3}{20G_1^3};
$$
then
\begin{equation}
\label{equation: theta nu}
\theta\nu=-3p_1\nu-\frac 32\lambda\nu-2\mu^2-p_3,
\end{equation}
where
$$
p_3=\frac{-10G_3'''+40G_3G_3''+21(G_3')^2-54G_3^2G_3'+9G_3^4}{50G_1^4}
$$
is a function invariant under~$\Gamma$ by Lemma~\ref{lemma: invariance of p}.
Using \eqref{equation: theta w}--\eqref{equation: theta nu} we find
\begin{align*}
\theta^2w
&=\biggl(\frac34\lambda^2-p_1\lambda-p_2+\mu\biggr)w,
\\
\theta^3w
&=\biggl(\frac38\lambda^3-\frac94p_1\lambda^2
+\biggl(p_1^2-\theta p_1-\frac52p_2+\frac32\mu\biggr)\lambda
\\ &\qquad
-\theta p_2+p_1p_2-3p_1\mu+\nu\biggr)w,
\displaybreak[2]\\
\theta^4w
&=\biggl(\frac3{32}\lambda^4-\frac94p_1\lambda^3
+\biggl(\frac{21}4p_1^2-3\theta p_1-3p_2+\frac34\mu\biggr)\lambda^2
+\dotsb\biggr)w,
\\ \intertext{and}
\theta^5w
&=\biggl(-\frac{15}{16}p_1\lambda^4
+\biggl(\frac{75}8p_1^2-\frac{15}4\theta
p_1-\frac{15}8p_2\biggr)\lambda^3+\dotsb\biggr)w.
\end{align*}
It follows that
\begin{align*}
\theta^5w
&=-10p_1\theta^4w
+\biggl(-\frac18(105p_1^2+30\theta p_1+15p_2)\lambda^3+\dotsb\biggr)w
\\
&=-10p_1\theta^4w-(35p_1^2+10\theta p_1+5p_2)\theta^3w
\\ &\qquad
+\biggl(-\frac18(300p_1^3+270p_1\theta p_1+30\theta^2p_1
+180p_1p_2+45\theta p_2)\lambda^2+\dotsb\biggr)w
\\
&=-10p_1\theta^4w-(35p_1^2+10\theta p_1+5p_2)\theta^3w
\\ &\qquad
-\biggl(50p_1^3+45p_1\theta p_1+5\theta^2p_1+30p_1p_2
+\frac{15}2\theta p_2\biggr)\theta^2w+\dotsb.
\end{align*}
Continuing this way, we find that $w$ satisfies the differential
equation given in the statement of Theorem~\ref{theorem: converse 2}.
\end{proof}

\begin{proof}[Proof of Theorems \rom{\ref{theorem: converse 3}} and
  \rom{\ref{theorem: converse 4}}]
Here we shall adopt all the notations in
the proof of Theorems \ref{theorem: converse 1} and~\ref{theorem: converse 2}.
In particular, we set
$$
  \lambda=\frac{G_2}{G_1}, \qquad
  \mu=\frac{4G_3'-3G_3^2}{20G_1^2}, \qquad
  \nu=\frac{4G_3''-10G_3'G_3+3G_3^3}{20G_1^3}.
$$
We first give a proof of Theorem~\ref{theorem: converse 4}.

Observe that
\begin{equation} \label{equation: w^2/G1}
  \begin{vmatrix}w_0&\theta w_0\\w_1&\theta w_1\end{vmatrix}
 =w_0^2z\frac{\d\tau_1}{\d z}=\frac{w_0^2}{G_1}.
\end{equation}
Setting $u=g\bigl|\begin{smallmatrix}w_0&\theta w_0\\w_1&\theta
  w_1\end{smallmatrix}\bigr|^{1/2}$, we have, by \eqref{equation: theta w},
\eqref{equation: p1 p2} and $\theta=z\d/\d z=G_1^{-1}\d/\d\tau$,
$$
  \theta u=u\left(\lambda-\frac{G_1'}{G_1^2}-2p_1\right)
 =u\left(\frac34\lambda-\frac14\rho+\frac52p_1\right),
$$
where we set
$$
  \rho=\frac{G_3}{G_1}.
$$
Then, by \eqref{equation: p1 p2},
\begin{equation} \label{equation: theta rho}
  \theta\rho=5\mu-p_1\rho-\frac12\lambda\rho+\frac14\rho^2.
\end{equation}
Using \eqref{equation: theta lambda}, \eqref{equation: theta mu},
\eqref{equation: theta nu} we find
\begin{align*}
\theta^2u
&=u\biggl(\frac38\lambda^2-\biggl(\frac14\rho+\frac92p_1\biggr)\lambda
-\frac14(10\theta p_1+3p_2-25p_1^2+2\mu-6p_1\rho\biggr),
\\
\theta^3u
&=u\biggl(\frac3{32}\lambda^3-\biggl(\frac3{32}\rho
+\frac{63}{16}p_1\biggr)\lambda^2
\\ &\qquad
-\frac3{16}\biggl(34\theta p_1+7p_2-109p_2^2+2\mu-14p_1\rho\biggr)\lambda+\cdots\biggr),
\\ \intertext{and}
\theta^4u
&=u\biggl(-\frac32p_1\lambda^3-\frac3{16}\biggl(38\theta p_1+5p_2-149p_1^2-8p_1\rho\biggr)\lambda^2
+\cdots\biggr).
\end{align*}
{}From these computations, we see that
$g\bigl|\begin{smallmatrix}w_0&\theta w_0\\w_1&\theta
w_1\end{smallmatrix}\bigr|^{1/2}$ satisfies the differential
equation in Theorem~\ref{theorem: converse 4}. By a similar argument,
we can show that the other three functions also satisfy the same
differential equation.

We now prove the monodromy part of Theorem \ref{theorem: converse 3}.
Set
$$
  u_0=\begin{vmatrix}w_0&\theta w_0\\w_1&\theta w_1\end{vmatrix}^{1/2}, \quad
  u_1=-\tau_2'u_0, \quad u_2=(\tau_3-\tau_2\tau_2')u_0, \quad
  u_3=(\tau_2-\tau_1\tau_2')u_0.
$$
According to Remark~\ref{r8}, up to a sign,
these functions are the same as the four
functions given in the statement of Theorem~\ref{theorem: converse 3}. Let
$$
 \gamma
 =\begin{pmatrix}A&B\\C&D\end{pmatrix}
 =\M{{a_{11}}{a_{12}}{a_{13}}{a_{14}}}
  {{a_{21}}{a_{22}}{a_{23}}{a_{24}}}
  {{a_{31}}{a_{32}}{a_{33}}{a_{34}}}
  {{a_{41}}{a_{42}}{a_{43}}{a_{44}}}
 \in\Gamma.
$$
It suffices to show that under the action of~$\gamma$,
\begin{equation} \label{equation: goal}
  \begin{pmatrix}u_3\\u_2\\u_1\\u_0\end{pmatrix}\mapsto
  \epsilon\gamma
  \begin{pmatrix}u_3\\u_2\\u_1\\u_0\end{pmatrix}
\end{equation}
for some complex scalar $\epsilon$ depending on~$\gamma$.

Using \eqref{equation: w^2/G1} and \eqref{equation: p 1} we find that
under the action of~$\gamma$,
$$
  u_0^2=\frac{w_0^2}{G_1}
  \mapsto{\chi(\gamma)^2(\fc t+\fd)^2}\frac{w_0^2}{G_1}.
$$
Recalling that $t=-\tau_2'$ by our formulation we have
$$
  \fc t+\fd
  =a_{41}(-\tau_1\tau_2+\tau_2)+a_{42}(-\tau_2\tau_2+\tau_3)
  -a_{43}\tau_2'+a_{44}.
$$
It follows that
$$
u_0\mapsto\epsilon(a_{41}u_3+a_{42}u_2+a_{43}u_1+a_{44}u_0)
$$
under the action of~$\gamma$ for some scalar~$\epsilon$.

For the behavior of $u_1$ under~$\gamma$,
we use~\eqref{equation: t}. We find
$$
  u_1=tu_0\mapsto\frac{\fa t+\fb}{\fc t+\fd}\left(\epsilon
  (\fc t+\fd)u_0\right)=\epsilon
  (a_{31}u_3+a_{32}u_2+a_{33}u_1+a_{34}u_0).
$$
For $u_2$ and $u_3$, we have
$$
  \begin{pmatrix}u_3\\ u_2\end{pmatrix}
 =\Tau\begin{pmatrix}t\\1\end{pmatrix}u_0.
$$
Under the action of $\gamma$ we have
$\Tau\mapsto(A\Tau+B)(C\Tau+D)^{-1}$ and
by~\eqref{equation: t}
$$
  \begin{pmatrix}t\\1\end{pmatrix}u_0\mapsto
  \epsilon\begin{pmatrix}\fa t+\fb\\\fc t+\fd\end{pmatrix}u_0
 =\epsilon(C\Tau+D)\begin{pmatrix}t\\ 1\end{pmatrix}u_0.
$$
It follows that
$$
  \begin{pmatrix}u_3\\u_2\end{pmatrix}\mapsto
  \epsilon(A\Tau+B)\begin{pmatrix}t\\1\end{pmatrix}u_0
$$
under the action of $\gamma$. From this we deduce that
$$
  u_2\mapsto\epsilon(a_{21}u_3+a_{22}u_2+a_{23}u_1+a_{24}u_0),
  \quad
  u_3\mapsto\epsilon(a_{11}u_3+a_{12}u_2+a_{13}u_1+a_{14}u_0).
$$
This establishes \eqref{equation: goal} and completes the proof of the theorems.
\end{proof}

\section{Guillera's generalization of Ramanujan's formulas for $1/\pi$}
\label{s5}

The modular parameterization of solutions of Picard--Fuchs linear
differential equations of order~3 has another curious application
to proving Ramanujan's series for $1/\pi$ \cite{Ra}, like
$$
\sum_{n=0}^\infty\frac{(4n)!}{n!^4}
(26390n+1103)\cdot\frac1{396^{4n}}
=\frac{99^2}{2\pi\sqrt2}.
$$
Note that the series on the left-hand side is a $\mathbb Q$-linear
combination of the ${}_3F_2$ hypergeometric series satisfying~\eqref{3F2}
and its derivative at a point close to the origin. The paper~\cite{Zu}
reviews ideas of proofs of Ramanujan's series and its several generalizations.

As already mentioned in Section~\ref{s1},
modular Picard--Fuchs differential equations of order~3
always come as the symmetric square of equations of order~2.
We consider such a second order differential equation
and proceed in the notation of Section~\ref{s1} until
equality~\eqref{e02}, where we choose $C=2\pi i$.

{}From \eqref{e02a} we have
$$
\frac{\d\tau}{\d z}
=\frac{W(u_0,u_1)}{u_0^2}
=\frac1{2\pi ig_0u_0^2},
$$
hence
\begin{equation}
\delta=\frac1{2\pi i}\,\frac{\d}{\d\tau}
=\frac1{2\pi i}\biggl(\frac{\d\tau}{\d z}\biggr)^{-1}\frac{\d}{\d z}
=g_0u_0^2\frac{\d}{\d z}.
\label{e03}
\end{equation}

Our next object is the function
\begin{equation}
v=v(\tau)=\delta\log u_0.
\label{e04}
\end{equation}
{}From \eqref{e03} we obtain
\begin{equation}
v=\frac{\delta u_0}{u_0}=g_0(z)u_0u_0'.
\label{e05}
\end{equation}

\begin{Lemma}
\label{c4}
The following functional equation is valid
for any $\gamma\in\Gamma$:
\begin{equation}
v(\gamma\tau)
=(c\tau+d)^2v(\tau)+\frac1{2\pi i}c(c\tau+d).
\label{e06}
\end{equation}
\end{Lemma}

\begin{demo}
Indeed,
\begin{align*}
v(\gamma\tau)
&=\frac1{2\pi i}\,\frac{\d}{\d(\gamma\tau)}\log(cu_1+du_0)
\\
&=\biggl(\frac{\d(\gamma\tau)}{\d\tau}\biggr)^{-1}
\cdot\frac1{2\pi i}\,\frac{\d}{\d\tau}
\log\bigl(u_0\cdot(c\tau+d)\bigr)
\\
&=(c\tau+d)^2\cdot\biggl(\delta\log u_0
+\frac1{2\pi i}\,\frac c{c\tau+d}\biggr)
\\
&=(c\tau+d)^2v(\tau)+\frac1{2\pi i}c(c\tau+d).
\end{align*}
\end{demo}

\begin{Lemma}
\label{c5}
For any integer $N\ge2$, the function
$\wt v(\tau)=v(\tau)-Nv(N\tau)$
is a $\Gamma'$-modular form of weight~$2$,
where
$$
\Gamma'=\Gamma_N'
=\biggl\{\gamma=\begin{pmatrix} a & b \\ cN & d \end{pmatrix}\in\Gamma:
\gamma^*=\begin{pmatrix} a & bN \\ c & d \end{pmatrix}\in\Gamma\biggr\}.
$$
\end{Lemma}

A consequence of this lemma is that
$\wt v=g_1(z)u_0^2$, where $g_1$ is an algebraic
function of~$z$.

\begin{demo}
For any
$\gamma=\bigl(\begin{smallmatrix} a & b \\ cN & d \end{smallmatrix}\bigr)\in\Gamma'$
we have
$$
N\cdot\gamma\tau=N\cdot\frac{a\tau+b}{cN\tau+d}
=\gamma^*(N\tau),
$$
hence, by~\eqref{e06},
\begin{align*}
\wt v(\gamma\tau)
&=v(\gamma\tau)-Nv(N\cdot\gamma\tau)
=v(\gamma\tau)-Nv(\gamma^*(N\tau))
\\
&=(cN\cdot\tau+d)^2v(\tau)+\frac1{2\pi i}cN(cN\cdot\tau+d)
\\ &\qquad
-N\biggl((c\cdot N\tau+d)^2v(N\tau)+\frac1{2\pi i}c(c\cdot N\tau+d)\biggr)
\\
&=(cN\tau+d)^2\bigl(v(\tau)-Nv(N\tau)\bigr)
=(cN\tau+d)^2\wt v(\tau)
\end{align*}
that implies the desired assertion.
\end{demo}

It is now time to glue the gathered information.
Take a quadratic irrationality $\tau_0$ (from the
upper half-plane) and an element $\gamma_0\in\Gamma$
such that
\begin{equation}
\gamma_0\tau_0=N\tau_0
\quad\text{for some integer $N\ge2$}.
\label{gamma-tau-N}
\end{equation}

\begin{Lemma}
\label{c6}
In the above notation,
the number $z_0=z(\tau_0)=z(\gamma_0\tau_0)$
is algebraic.
\end{Lemma}

\begin{demo}
This follows from the fact that $z(\tau)$ and
$z(N\tau)$ are connected by a (modular) polynomial
equation with integer coefficients.
Substituting $\tau=\tau_0$ into this equation
gives a polynomial for the number $z(\tau_0)=z(N\tau_0)$.
\end{demo}

{}From
$v(\tau_0)-Nv(N\tau_0)=\wt v(\tau_0)
=g_1(z_0)u_0(z_0)^2$ and Lemma~\ref{c4} it follows that
\begin{align*}
v(N\tau_0)
&=v(\gamma_0\tau_0)
=(c_0\tau_0+d_0)^2v(\tau_0)+\frac1{2\pi i}c_0(c_0\tau_0+d_0)
\\
&=(c_0\tau_0+d_0)^2\bigl(Nv(N\tau_0)+g_1(z_0)u_0(z_0)^2\bigr)
+\frac1{2\pi i}c_0(c_0\tau_0+d_0),
\end{align*}
hence
\begin{equation}
v(N\tau_0)
=\frac{(c_0\tau_0+d_0)^2g_1(z_0)u_0(z_0)^2+c_0(c_0\tau_0+d_0)/2\pi i}
{1-N\cdot(c_0\tau_0+d_0)^2}.
\label{e07}
\end{equation}
On the other hand, from~\eqref{e05} we have
\begin{equation}
v(N\tau_0)
=\frac12g_0(z)\bigl(u_0(z)^2\bigr)'\big|_{z=z_0}.
\label{e08}
\end{equation}
It remains to eliminate $v(N\tau_0)$ in~\eqref{e07}
and~\eqref{e08}:
\begin{equation}
\begin{split}
\frac1\pi
&=2i\frac{1-N\cdot(c_0\tau_0+d_0)^2}{2c_0(c_0\tau_0+d_0)}
g_0(z_0)\bigl(u_0(z)^2\bigr)'\big|_{z=z_0}
\\ &\qquad
-2i\frac{c_0\tau_0+d_0}{c_0}g_1(z_0)\bigl(u_0(z)^2\bigr)\big|_{z=z_0}.
\end{split}
\label{e09}
\end{equation}
This is a Ramanujan-type series for $1/\pi$.

\begin{Remark}
\label{r11}
There are two places, where the use of modularity is crucial:
the algebraicity of $z_0=z(\tau_0)$ (Lemma~\ref{c6}) and
the algebraicity of $g_1(z)=\wt v/u_0^2$ (Lemma~\ref{c5}).
The fact that $\tau$~is algebraic (and even quadratic)
follows from $\gamma_0\tau_0=N\tau_0$, while the algebraicity
of $g_0(z)$ is a purely analytic fact (see the analytic proof
of Corollary in Section~\ref{s1}). It would be nice to avoid the modularity
completely, thus providing a purely differential equation
proof of equality~\eqref{e09}.
\end{Remark}

It seems that there exist analogous algebraic relations
in the case of Picard--Fuchs fourth and fifth order linear
differential equations considered above.

The differential equation~\eqref{e23} in Example~\ref{ex3}
and its analytic solution $w_0=w_0(z)$, which is a ${}_5F_4$ hypergeometric
series, are related to the following formulas for $1/\pi^2$
proved recently by J.~Guillera~\cite{G1},~\cite{G2}:
\begin{align}
\refstepcounter{equation}
\sum_{n=0}^\infty\binom{2n}n^5
(20n^2+8n+1)\biggl(-\frac1{2^{12}}\biggr)^n&=\frac8{\pi^2},
\tag{\theequation a}
\label{JGa}
\\
\sum_{n=0}^\infty\binom{2n}n^5
(820n^2+180n+1)\biggl(-\frac1{2^{20}}\biggr)^n&=\frac{128}{\pi^2}.
\tag{\theequation b}
\label{JGb}
\end{align}
Namely, following the notations in Example \ref{ex3}, for the two
specializations of~$z$,
\begin{itemize}
\item[(a)] $z=-1/2^{12}$, and
\item[(b)] $z=-1/2^{20}$,
\end{itemize}
we discovered experimentally that
\begin{gather}
\refstepcounter{equation}
3\tau_1+4\tau_2+4\tau_3-2(\tau_1\tau_3-\tau_2^2)=14,
\tag{\theequation a}
\label{e24a}
\\
\refstepcounter{equation}
(\tau_1-2)\frac{\d\tau_2}{\d\tau_1}-\tau_2=\sqrt 5+1,
\tag{\theequation a}
\label{e25a}
\displaybreak[2]\\ \intertext{and}
\addtocounter{equation}{-1}
7\tau_1+12\tau_2+4\tau_3-2(\tau_1\tau_3-\tau_2^2)=78,
\tag{\theequation b}
\label{e24b}
\\
\refstepcounter{equation}
(\tau_1-2)\frac{\d\tau_2}{\d\tau_1}-\tau_2=\sqrt{41}+3,
\tag{\theequation b}
\label{e25b}
\end{gather}
respectively, where $\tau_i=\tau_i(z)$ are defined as in Example~\ref{ex3}.
(One may also verify that, for an arbitrary $z<0$, we have
non-holomorphic relations $\Re(\tau_1/2-1)=0$, $\Im(\tau_1/2+\tau_2)=0$,
$\Re(\tau_1/4+\tau_2+\tau_3-1)=0$, and $\Re(\d\tau_2/\d\tau_1+1/2)=0$.)
Equations~\eqref{e24a} and~\eqref{e24b} show that the points
$\Tau=\left(\begin{smallmatrix}\tau_1&\tau_2\\\tau_2&\tau_3\end{smallmatrix}\right)$
lie on certain rational Humbert surfaces;
the equations can be written as
\begin{align}
\refstepcounter{equation}
\det(\gamma^{\text{(a)}}\Tau)&=-5\phantom0
\quad\text{or}\quad
\gamma_1(\gamma^{\text{(a)}}\Tau)=\frac15\gamma^{\text{(a)}}\Tau,
\tag{\theequation a}
\label{e26a}
\\
\det(\gamma^{\text{(b)}}\Tau)&=-41
\quad\text{or}\quad
\gamma_1(\gamma^{\text{(b)}}\Tau)=\frac1{41}\gamma^{\text{(b)}}\Tau,
\tag{\theequation b}
\label{e26b}
\end{align}
respectively, where
\begin{equation*}
\gamma^{\text{(a)}}=\begin{pmatrix} 1 & 0 & -2 & 1 \\ \frac12 & 1 & 0 & -1 \\
0 & 0 & 1 & -\frac12 \\ 0 & 0 & 0 & 1 \end{pmatrix}
\quad\text{and}\quad
\gamma^{\text{(b)}}=\begin{pmatrix} 1 & 0 & -2 & 3 \\ \frac12 & 1 & 2 & -2 \\
0 & 0 & 1 & -\frac12 \\ 0 & 0 & 0 & 1 \end{pmatrix}
\end{equation*}
are matrices in $\Sp_4(\mathbb Q)$, and $\gamma_1$~is defined in~\eqref{gamma01}.
Relations~\eqref{e26a} and~\eqref{e26b} may be viewed as an analogue of~\eqref{gamma-tau-N},
although the matrices $\gamma^{\text{(a)}}$ and $\gamma^{\text{(b)}}$ do not
belong to the monodromy group of Example~\ref{ex3}.

The equalities~\eqref{e24a} and~\eqref{e24b} happen to hold for the non-holomorphic
embedding $\Zeta$ in Section~\ref{s-klemm} as well: one simply replaces
the corresponding entries of~$\Tau$ by~$\Zeta$.

We have verified the five other examples in~\cite{G2} and conclude that
the algebraicity seems to appear in all Guillera's identities for $1/\pi^2$
(both conjectural and proved).
What is a theoretic background for these algebraic relations?

\medskip
\noindent
\textbf{Acknowledgments.}
We thank Gert Almkvist, Daniel Bertrand, Frits Beukers, J\'esus Guillera,
Albrecht Klemm, Duco van Straten, Noriko Yui, and Don Zagier
for useful conversations on parts of this work.
The authors are thankful to the staff
of the Max Planck Institute for Mathematics in Bonn
for the wonderful working conditions they experienced doing this research.

\vskip7mm

\section*
{\textbf{Appendix. \uppercase{On a subgroup of infinite index in $\Sp_4(\Z)$}}}
\label{sA}

\medskip
\hbox to\hsize{\hss \uppercase{Vicen\c tiu Pa\c sol}\hss}

\makeatletter

\let\@makefnmark\relax  \let\@thefnmark\relax
\@footnotetext{\def\par{\let\par\@par}%
  The work was supported by a fellowship
  of the Max Planck Institute for Mathematics (Bonn).}

\makeatother

\vskip5mm

Let $T=\bigl(\begin{smallmatrix} 1&1\\0&1\end{smallmatrix}\bigr)$ and
$S=\bigl(\begin{smallmatrix} 0&1\\-1&0\end{smallmatrix}\bigr)$ be the
two generators for $\SL_2(\Z)$.

Let $\Gamma$ be the subgroup of $\Sp_4(\Z)$ generated by the
matrices
$$
\gamma_0=\begin{pmatrix}-STS&M\\0&T\end{pmatrix}
\quad\text{and}\quad
\gamma_1=\begin{pmatrix}0&S\\-S&0\end{pmatrix},
$$
where $M=\bigl(\begin{smallmatrix}4&2\\-2&1\end{smallmatrix}\bigr)$.
These are exactly the matrices in~\eqref{gamma01}.

\begin{ATheorem}
\label{th:inf}
The group $\Gamma$ has infinite index in $\Sp_4(\Z)$.
\end{ATheorem}

The idea of the proof is to find a principal $\Sp_4(\Z)$-module such
that $\Gamma$ has infinitely many orbits.

\medskip
For a vector $\ba=(a,b,c,d)\in\Z^4$, we say it is reduced if it is
primitive (i.e., $\gcd(a,b,c,d)=1$) and it satisfies the following
conditions:
\begin{enumerate}
\item[(1)] $a\ge 0$, $b\le 0$, $c\ge 0$, $d\ge 0$,
\item[(2)] $-b\le d/2-c$, and
\item[(3)] $c\le a/2+b$.
\end{enumerate}
In particular, we have
\begin{enumerate}
\item[(4)] $-b\le a/2$ and
\item[(5)] $c\le d/2$.
\end{enumerate}

On the set of primitive vectors $(\Z^4)'$ we introduce the
involution $\varepsilon$ by the rule
$$
\varepsilon(a,b,c,d):=(d,-c,-b,a).
$$
One easily check that $\varepsilon$ preserves the reduced vectors.

Finally, we introduce an algorithm that produces reduced vectors
starting with any primitive one.

\begin{algo}
Fix a pair $\ba=(a,b,c,d)\in(\Z^4)'$. Consider the following
process.

\noindent{\sl Step\/}~0.
Put $\ba_0=(a_0,b_0,c_0,d_0):=(|a|,-|b|,|c|,|d|)$.

\noindent{\sl Step\/}~1.
Let $b_1:=-\min_{n\in\Z}|b_0+na_0|$. Put $\ba_1:=(a_0,b_1,c_0,d_0)$.

\noindent{\sl Step\/}~2.
Let $c_1:=\min_{n\in\Z}|c_0+nd_0|$. Put $\ba_2:=(a_0,b_1,c_1,d_0)$.

\noindent{\sl Step\/}~3.
Let $n_0\in\Z$ such that $\min_{n\in\Z}|b_1+n(d_0+2c_1)|=|b_1+n_0(d_0+2c_1)|$.
Put $b_2:=-|b_1+n_0(d_0+2c_1)|$, $a_2=|a_0+n_0(4c_1-2d_0)|$, and
$\ba_3:=(a_2,b_2,c_1,d_0)$.

\noindent{\sl Step\/}~4.
Let $n_1\in\Z$ such that $\min_{n\in\Z}|c_1+n(-a_2+2b_2)|=|c_1+n_1(a_2+2b_2)|$.
Put $c_2:=|c_1+n_1(-a_2+2b_2)|$, $d_2=|d_0-n_1(2a_2+4b_2)|$, and
$\ba_4:=(a_2,b_2,c_2,d_2)$.

\noindent{\sl Step\/}~5.
If $\ba_4\ne\ba$ repeat from Step~1 starting with $\ba_4$.

The process must ends by the well order principle. We call the end
result by $r(\ba)$.
\end{algo}

\begin{ALemma}
\label{lA1}
We have the following properties for this process:
\begin{enumerate}
\item $r(\varepsilon(\ba))=\varepsilon(r(\ba))$,
\item $r(a,b,c,d)=r(\pm a,\pm b,\pm c,\pm d)=r(r(a,b,c,d))$,
\item $\ba$ is reduced if and only if $r(\ba)=\ba$.
\end{enumerate}
\end{ALemma}

\begin{proof}
The first two properties are obvious from the definitions.

For the third one, if $\ba$ is reduced, none of Steps~0--4
of the algorithm produce anything new (easy check).
Therefore, the process ends and $r(\ba)=\ba$.

Conversely, if $r(\ba)=\ba$, we must have that none of Steps~0--4
of the algorithm produce anything new, since the
algorithm decreases (not necessarily strict) the absolute value of
the components at each step. Conditions (1)--(5) in the
definition of a reduced vector are easily checked to be satisfied.
\end{proof}

Lastly, we need a new definition.
For two primitive vectors $\ba_1,\ba_2$ we say they are equivalent,
$\ba_1\sim\ba_2$, if $r(\ba_1)=r(\ba_2)$ or
$r(\ba_1)=\varepsilon(r(\ba_2))$.
Notice that this relation is an equivalence relation by Lemma~\ref{lA1}.

\medskip
We give now the main result which implies Theorem~\ref{th:inf}.

We let $\GL_4(\Z)$ act on $(\Z^4)'$ by multiplication on
the left, so the vectors are considered as column vectors.

Let $\wt\Gamma$ be the subgroup of $\GL_4(\Z)$ generated by
the following matrices:
\begin{gather*}
A=\begin{pmatrix}-STS&0\\0&E\end{pmatrix}
=\begin{pmatrix}
1 & 0 & 0 & 0 \\ -1 & 1 & 0 & 0 \\
0 & 0 & 1 & 0 \\ 0 & 0 & 0 & 1
\end{pmatrix}, \quad
B=\begin{pmatrix}E&0\\0&T\end{pmatrix}
=\begin{pmatrix}
1 & 0 & 0 & 0 \\ 0 & 1 & 0 & 0 \\
0 & 0 & 1 & 1 \\ 0 & 0 & 0 & 1
\end{pmatrix},
\\
C=\begin{pmatrix}E&-ST^{-1}SMT^{-1}\\0&E\end{pmatrix}
=\begin{pmatrix}
1 & 0 & 4 & -2 \\ 0 & 1 & 2 & 1 \\
0 & 0 & 1 & 0 \\ 0 & 0 & 0 & 1
\end{pmatrix},
\quad\text{and}\quad \gamma_1,
\end{gather*}
where $E$ is the $2\times2$ identity matrix.
Since $ACB=\gamma_0$, we do have $\Gamma\subseteq\wt\Gamma$.

\begin{ALemma}
\label{lA2}
We have $\ba\sim \gamma\ba$ for all $\gamma\in \wt\Gamma$.
\end{ALemma}

\begin{proof}
It is obvious that it is enough to prove this only for the
generators of~$\wt\Gamma$.

If $\gamma=\gamma_1$ we have
$$
r(S\ba)=r(\varepsilon(\ba)))=\varepsilon(r(\ba)).
$$

If $\gamma=A$ we have
$$
r\left(A\begin{pmatrix}a\\b\\c\\d\end{pmatrix}\right)
=r\left(\begin{pmatrix}a\\b-a\\c\\d\end{pmatrix}\right)
=r(\ba).
$$
The last equality follows from the fact that after Step~1
both quantities are equal, i.e.,
$$
\min_{n\in\Z}|b+na|=\min_{n\in\Z}|na\pm(b-a)|.
$$

If $\gamma=B$, same but for $c$ and $d$ reasoning.

If $\gamma=C$, Step~3 takes care of the equality.
\end{proof}

\begin{ALemma}
\label{lA3}
If $\ba_1,\ba_2$ are reduced and $\ba_1\ne\ba_2$ and
$\ba_1\ne\varepsilon(\ba_2)$, then $\ba_1\not\sim_\Gamma\ba_2$.
\end{ALemma}

\begin{proof}
The fact the vectors are reduced implies that they are equal to
their reduced vectors. The condition in the corollary implies that
$\ba_1\not\sim\ba_2$, therefore, by the previous result we cannot
have equivalence under $\wt\Gamma$, in particular under $\Gamma$.
\end{proof}

To end the proof of our main result, we just have to observe that
we have infinite number of vectors which are reduced and
nonequivalent. For example, take the vectors $\ba_{p,q}:=(p,-1,0,q)$
for $p,q\ge2$ any two positive integers.

In this case, since $\Sp_4(\Z)$ acts transitively on the set of
primitive vectors in $\Z^4$ (so $(\Z^4)'$ is the principal
$\Sp_4(\Z)$-module we are considering), let $\gamma_{p,q}\in\Sp_4(\Z)$
such that $\gamma_{p,q}\cdot\trans{(1,0,0,0)}=\ba_{p,q}$. Then
$\gamma_{p,q}\not\sim_\Gamma\gamma_{p',q'}$ for any $(p,q)\ne (p',q')$.
Otherwise, $\ba_{p',q'}=\gamma\cdot\ba_{p,q}$ for some $\gamma\in\Gamma$. By
Lemma~\ref{lA3} this is impossible.

\medskip
{\footnotesize
\textsc{Max-Planck-Institut f\"ur Mathematik, Vivatsgasse 7, 53111 Bonn, Germany}

\textit{E-mail address}: \href{mailto:vpasol@gmail.com}{\texttt{vpasol@gmail.com}}
}


\end{document}